\documentclass[a4paper,12pt,french]{article}
\usepackage[latin1]{inputenc}
\usepackage{babel}
\usepackage[dvips]{graphicx}
\usepackage{epsfig}
\usepackage[usenames,dvipsnames]{color}
\usepackage{amssymb}
\title {MATRICE MAGIQUE ASSOCI\'EE \`A UN GERME DE COURBE PLANE ET DIVISION PAR L'ID\'EAL JACOBIEN}

\newcommand{\stl}{\vspace{3mm}}

\newtheorem{defi}{\sc Definition}[section]
\newtheorem{theo}[defi]{\sc Th\'eor\`eme}
\newtheorem{prop}[defi]{\sc Proposition}
\newtheorem{nota}[defi]{\sc Notation}
\newtheorem{coro}[defi]{\sc Corollaire}
\newtheorem{rema}[defi]{\sc Remarque}
\newtheorem{lemm}[defi]{\sc Lemme}
\newtheorem{exem}[defi]{\sc Exemple}
\newtheorem{form}[defi]{\sc Formulaire}

\newenvironment{demo}[1][]{\noindent {\it Preuve.} \protect\nopagebreak
       \rm #1}{\protect\nopagebreak $\square $ \par\stl}

\begin{document}
\begin{center}
 { \Large \bf MATRICE MAGIQUE ASSOCI\'EE 

\`A UN GERME 
DE COURBE PLANE  ET \\

DIVISION PAR L'ID\'EAL JACOBIEN}
\end{center}

\stl

\begin{center}
  Joël Briançon, \  Philippe Maisonobe \&
  Tristan Torrelli\footnote{
\{briancon, phm, torrelli\}@math.unice.fr, 
Laboratoire 
J.A. Dieudonn\'e,
U.M.R. du C.N.R.S. 6621,
Universit\'e de Nice Sophia-Antipolis,
Parc Valrose, F-06108 Nice Cedex~2.

\noindent Classification ${\cal AMS}$ 2000 : 32S40, 32S10, 32C38, 32C40,
14B05.

\noindent Mots clef : Polyn\^ome de Bernstein-Sato, 
germe de courbe plane, matrices magiques.}
\end{center}

\stl

\noindent{\sc R\'esum\'e}. Dans l'anneau des
germes de fonctions holomorphes \`a l'origine de 
${\bf C}^2$, nous nous donnons une fonction $f$
d\'efinissant une singularit\'e isol\'ee, et nous nous
int\'eressons \`a l'\'equation : $uf'_x+vf'_y=wf$, 
lorsque la fonction  $w$ est donn\'ee. Nous introduisons les
multiplicit\'es d'intersection relatives de $w$ et
$f'_y$ le long des branches de $f$ et nous \'etudions
les solutions \`a l'aide de ces valuations. Gr\^ace
aux r\'esultats ainsi d\'emontr\'es, nous construisons
explicitement une \'equation fonctionnelle v\'erifi\'ee
par $f$.

\stl

\ 

\smallskip 

\noindent{\large \bf Introduction}

\stl

Soit $\overline{K}=\bigcup_{d\in{\bf N}^*} 
{\bf C}[[x^{1/d}]][1/x]$ le corps des s\'eries de Puiseux
muni de sa valuation naturelle $\nu$, et 
$f=\prod_{i=1}^n(y-a_i)$ un polyn\^ome unitaire,
r\'eduit, de degr\'e $n\geq 2$, \`a coefficient dans
$\overline{K}$. Nous lui associons la famille des
multiplicit\'es~:
$m_{i,j}=\nu(a_i-a_j)$ pour $1\leq i\not= j\leq n$,
et la matrice magique  $A=(\alpha_{i,j})$
d\'efinie par : $\alpha_{i,j}=-m_{i,j}$ pour $i\not=j$
et $\alpha_{i,i}=\sum_{j\not=i}m_{i,j}$. Dans un
premier temps, nous d\'emontrons que cette
matrice est diagonalisable sur ${\bf Q}$ et nous
exhibons ses valeurs propres et ses sous-espaces 
propres (Corollaire \ref{AFcourbliss}).

\stl

Nous nous proposons alors d'expliciter une solution
$(u,v)$ du syst\`eme lin\'eaire $uf'_x+vf'_y=wf$
lorsque $w$ est un polyn\^ome donn\'e de
$\overline{K}[y]$, de degr\'e strictement inf\'erieur
\`a $n$ ; notons $\overline{\cal E}=\overline{K}[y]^{(n)}$
cet espace de polyn\^omes, et 
$\varepsilon_i=\prod_{j\not= i}(y-a_j)$, $1\leq i\leq n$,
la base d'interpolation de Lagrange (aux coefficients pr\`es).
Nous filtrons $\overline{\cal E}$ \`a l'aide de la valuation des
coordonn\'ees dans cette base~:
$$\mbox{val}(w)=
\mbox{inf}\left\{\nu\left(\frac{w(x,a_i)}{f'_y(x,a_i)}\right)\ ;\ 1\leq i\leq n\right\} \, \mbox{pour }w\in\overline{\cal E}\ .$$
Nous nous apercevons que les matrices colonnes des 
coordonn\'ees de $u$ et $w$ respectivement, dans cette
base, sont li\'ees par~:
$${\cal W}={\cal A}\cdot {\cal U}$$
o\`u $\cal A$ est une matrice magique \`a coefficients
dans $\overline{K}$, dont la forme initiale est $(1/x)A$.
Nous pouvons alors 
d\'eterminer la valuation de la
solution $(u,v)$ et sa forme initiale, sous
certaines conditions (Corollaire \ref{laresolwf}). 

\stl

Nous appliquons ces r\'esultats au cas o\`u 
$f\in{\bf C}\{x \}[y]$ est un polyn\^ome distingu\'e
de degr\'e $n$, d\'efinissant un germe de courbe
plane \`a singularit\'e isol\'ee ; nous en d\'eduisons
que lorsque $w\in{\bf C}\{x\}[y]^{(n)}$ est de valuation
positive ou nulle, $wf$ appartient \`a l'id\'eal jacobien
(Corollaire \ref{cororesolwf}). Nous retrouvons, par exemple, 
l'appartenance de $f^2$ \`a cet id\'eal. \'Egalement, nous 
g\'en\'eralisons ces r\'esultats au cas o\`u $f$ n'est plus 
un polyn\^ome r\'eduit.

\stl

En application, nous consacrons la derni\`ere partie 
 \`a la construction d'un multiple
du polyn\^ome de Bernstein de $f$ (\`a singularit\'e
isol\'ee de multiplicit\'e $n$) et de l'op\'erateur
associ\'e, v\'erifiant l'\'equation fonctionnelle~:
$$b(s)f^s=P\cdot f^{s+1}\ .$$
Le polyn\^ome trouv\'e $b(s)$ se calcule \`a 
partir des valeurs propres de la matrice magique 
de $f$, et il ne d\'epend donc que du type topologique
du germe de courbe d\'efini par $f$ (Th\'eor\`eme 
\ref{theomultopo}). Cela r\'epond \`a la question de
la construction effective, question qui interpelle 
les deux premiers auteurs depuis longtemps, apr\`es
la d\'emonstration de l'existence du polyn\^ome de
Bernstein par M. Kashiwara (\cite{K}) \`a l'aide de la
r\'esolution des singularit\'es.
 
\stl

\section{Matrices magiques} 

\subsection{La matrice magique associ\'ee \`a un bouquet}           
  
\subsubsection{Les bouquets}

Consid\'erons une famille de $n$ germes de courbes lisses 
distinctes ($n\geq 2$), transverses \`a l'axe des $y$ dans 
le plan ${\bf C}^{2}$ des couples $(x,y)$. Ce 
`bouquet de courbes lisses' est d\'efini par un unique polynôme 
unitaire $f$ appartenant \`a ${\bf C}\{x\}[y]$ :
$$f= \prod_{i=1}^n (y - a_i)$$
avec $a_i\in {\bf C}\{x\}$, $1\leq i\leq n$, deux-\`a-deux distincts.\\

Dans nos r\'esultats, les probl\`emes de convergence ne posent pas 
de difficult\'e s\'erieuse et nous prendrons les $a_{i}$ dans ${\bf C}[[x]]$. 
Plus g\'en\'eralement, nous envisagerons les $a_{i}$ dans le corps des 
fractions de ${\bf C}[[x]]$~:
$$K = {\bf C}((x)) = {\bf C}[[x]][1/x]$$
et dans la clôture alg\'ebrique de $K$~:
$$\overline{K}=\bigcup_{d\in {\bf N}^*} {\bf C}[[x^{1/d}]] [1/x]$$
que nous munissons de sa valuation naturelle 
$\nu: \overline{K} \rightarrow {\bf Q} \cup \{+\infty\}$. 

\begin{defi} \label{bouqpass}
  Un {\em bouquet de $n$ branches} est une famille ordonn\'ee de $n$
  \'el\'ements ($n\geq 2$) $a_i\in \overline{K}$, $1\leq i\leq n$, 
deux-\`a-deux distincts. 
  Le bouquet {\em passe par l'origine} lorsque les valuations $\nu(a_i)$,
$1\leq i\leq n$, sont strictement positives. 
\end{defi}

Un bouquet est donc donn\'e par un polyn\^ome unitaire
r\'eduit de degr\'e $n\geq 2$ de $\overline{K}[y]$, et avec un
ordre sur ses racines ; nous continuerons \`a \'ecrire :
$f=\prod_{i=1}^n(y-a_i)$.

\begin{exem}{\em \label{expoldist}
Soit $f\in{\bf C}[[x]][y]$ un polyn\^ome distingu\'e de
degr\'e $n$, c'est-\`a-dire unitaire et v\'erifiant : $f(0,y)=y^n$.
Notons $d\in {\bf N}^*$ le p.p.c.m des degr\'es de ses facteurs irr\'eductibles.
D'apr\`es le th\'eor\`eme de Newton-Puiseux, nous savons que~:
$f(t^d , y)=\prod_{i=1}^n (y - \phi_i(t))$ avec $\phi_i \in {\bf C}[[t]]$ 
et $\phi_i(0)=0$ pour tout $1\leq i\leq n$. Lorsque $f$ est r\'eduit, 
en num\'erotant ainsi les racines, nous obtenons le bouquet 
 passant 
par l'origine : $f= \prod_{i=1}^n (y - \phi_i(x^{1/d}))$.}
\end{exem}

Nous nous int\'eresserons aussi \`a des bouquets de branches
multiples.

\begin{defi} \label{defbouqmult}
     Un {\em bouquet de branches avec multiplicit\'es} est la
donn\'ee d'un bouquet de $n$-branches et d'une
famille $\mu=(\mu_1,\ldots,\mu_n)$ d'entiers naturels 
non nuls. 
\end{defi}

Nous \'ecrirons~: $f_\mu=\prod_{i=1}^n(y-a_i)^{\mu_i}$.

\stl

  Lorsque les $a_i$ appartiennent \`a $K$, 
nous parlerons alors de `bouquet de courbes m\'eromorphes'. 
\`A un bouquet donn\'e $f=\prod_{i=1}^n(y-a_i)\in \overline{K}[y]$,
nous associons  la famille
des multiplicit\'es d'intersection de ses branches deux-\`a-deux : 
$\{m_{i,j} = \nu (a_{i} - a_{j})\ ;\  i \not = j\}$
et nous posons : 
$m_{i} = \sum_{j \not=i} m_{i,j}=\nu(f'_y(a_{i}))$. 
En prenant les valuations dans la somme :
$$(a_{i} - a_{j}) + (a_{j} - a_{k}) + (a_{k} - a_{i}) = 0$$
pour des indices distincts $i,j,k$, nous constatons encore
que le minimum des entiers $ m_{i,j}, m_{j,k}, m_{k,i} $
est atteint pour au moins deux couples d'indices.

\subsubsection{La matrice magique associ\'ee \`a un bouquet 
de branches}

Rappelons qu'une matrice carr\'ee $A$ \`a coefficients dans un anneau
commutatif $\cal K$ est dite {\it magique} si la somme des coefficients 
de ses lignes et de ses colonnes est la même ; nous la notons $s(A)$. L'ensemble des matrices magiques 
$n\times n$ \`a coefficients dans $\cal K$ est une $\cal K$-alg\`ebre,
not\'ee ${\cal M}ag(n,\cal K)$,
et $$s: {\cal M}ag(n,\cal K) \longrightarrow \cal K $$
est un morphisme de $\cal K$-alg\`ebres. 
Le noyau de $s$, c'est-\`a-dire l'ensemble des matrices magiques de 
somme nulle, est un id\'eal bilat\`ere de
 ${\cal M}ag(n,\cal K)$ not\'e ${\cal M}ag_0(n,\cal K)$.\\

Lorsque $\cal K$ est totalement ordonn\'e, nous d\'efinissons  
l'ensemble ${\cal M}ag^{\star}(n,\cal K)$ des {\em matrices
magiques exceptionnelles} :  c'est l'ensemble des matrices magiques
$A=(\alpha_{i,j})$ v\'erifiant la propri\'et\'e : 

\stl

\noindent $(\star)$ : pour tout triplet $(i,j,k)$ d'indices distincts, 
le maximum de $\alpha_{i,j}, \alpha_{j,k},\alpha_{k,i}$ 
$\mbox{\ \ \ \ \ }$ est atteint au moins deux fois.\\
 
Signalons que cet ensemble ${\cal M}ag^{\star}(n,\cal K)$ n'est en
g\'en\'eral stable ni par  addition, ni par multiplication 
par un \'el\'ement de $\cal K$. Nous noterons enfin 
${\cal M}ag^{\star}_0(n,\cal K)$ l'ensemble 
des matrices exceptionnelles de somme nulle. 

\stl

Donnons maintenant la d\'efinition motivant ces rappels et ces 
notations. 

\begin{defi}
  Soit $f=\prod^n_{i=1}(y-a_i)\in\overline{K}[y]$ un polyn\^ome
r\'eduit de degr\'e $n\geq 2$. 
  La {\em matrice magique associ\'ee au bouquet} d\'efini par $f$ 
est la matrice carr\'ee $A=(\alpha_{i,j})$ \`a $n$ lignes 
et $n$ colonnes, \`a coefficients rationnels, d\'efinie par : 
$\alpha_{i,j} = - m_{i,j}$ pour $i\not=j$, et 
 $\alpha_{i,i} = m_{i} = \sum_{j \not=i} m_{i,j}$ pour $1\leq i\leq n$.
\end{defi}

 Ainsi, la matrice magique associ\'ee \`a un bouquet de courbes lisses est 
dans ${\cal M}ag^{\star}_0(n,\bf Z)$
 et dans ${\cal M}ag^{\star}_0(n,\bf Q)$ pour un bouquet de branches ; elle est de plus sym\'etrique.

\begin{rema}{\em
Lors de la permutation des branches d'un bouquet, la matrice
magique est remplac\'ee par la matrice semblable donn\'ee
par l'action de la matrice de permutation.}
\end{rema}

\begin{exem}{\em \label{exrec}
   La matrice magique associ\'ee au bouquet de courbes lisses d\'efini par : 
$$ f = y(y-1)(y-1-x)(y+1)(y+1+x)(y+1+x+x^{3}) $$ 
(avec $a_1=0$, $a_2=1$, ...) est : 

$$ A \; = \; \left(\begin{array}{rrrrrr}
0 & 0 & 0 & 0 & 0 & 0 \\
0 & 1 &-1 & 0 & 0 & 0 \\
0 &-1 & 1 & 0 & 0 & 0 \\
0 & 0 & 0 & 2 &-1 &-1 \\
0 & 0 & 0 &-1 & 4 &-3 \\
0 & 0 & 0 &-1 &-3 & 4
\end{array} \right) $$}
\end{exem}

Dans le cas d'un bouquet de branches avec multiplicit\'es, 
nous consid\`ererons une matrice $A_\mu$ d\'efinie
\`a partir des multiplicit\'es $m_{i,j}$ et $\mu$, et
qui  g\'en\'eralise la matrice magique $A$. Le paragraphe
\ref{magikgen} lui est enti\`erement consacr\'e.

\subsection{L'arbre associ\'e \`a un bouquet}

\subsubsection{Le cas d'un bouquet de courbes lisses}

Rappelons que la donn\'ee des multiplicit\'es d'intersection 
$\{m_{i,j} \ ;\  i\neq j\}$ d'un bouquet de $n$ courbes lisses 
\'equivaut \`a la donn\'ee de l'arbre des points infiniment 
voisins (ou `arbre d'\'eclatements', ou `arbre de d\'esingularisation'). 

\,

\stl

\begin{figure}[h]
\begin{center}
\includegraphics[height=6cm]{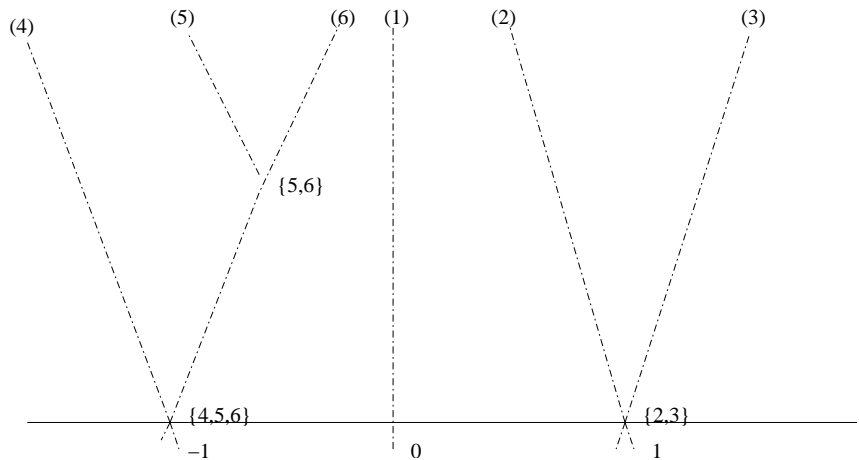}

\end{center}
\caption{L'arbre d'\'eclatements de l'exemple \ref{exrec}}
\end{figure}

 Une {\em bifurcation} est un point infiniment voisin dont l'\'eclatements
 s\'epare des branches du bouquet ; en particulier, elle est caract\'eris\'ee 
par les branches qui y passent. Nous allons d\'ecrire cet arbre 
d'\'eclatement en pr\'ecisant pour chaque bifurcation les indices 
des branches qui y passent et sa hauteur.

\stl

 Soit $\cal B $ l'ensemble des parties $T$  
de $\{1,2,\ldots,n\}$ v\'erifiant les propri\'et\'es suivantes :\\
\indent - le nombre d'\'el\'ements de $T$ est au moins \'egal \`a deux ;\\
\indent - pour tout couple d'indices distincts $(i,j)\in T^2$ : 
$m_{i,j}>0$ ;\\
\indent  - pour tout couple d'indices distincts $(i,j)\in T^2$
 et pour tout indice $k$ n'appartenant pas \`a $T$ :
  $ m_{i,j}$ est strictement sup\'erieur \`a $ m_{i,k} = m_{j,k}$.\\

\`A une bifurcation de l'arbre d'\'eclatement, on associe alors l'ensemble 
$T\subset{\{1,2,\ldots, n\}}$ des indices des branches qui y passent ; 
et la bifurcation correspondant \`a $T \in{\cal B}$ est le point infiniment voisin dont l'\'eclatement s\'epare des branches index\'ees par $T$.\\

Pour $T \in{\cal B}$, nous notons~:
$$\alpha(T) = \mbox{inf}\{ m_{i,j} \ ; \ i\in T , j\in T, i\neq j \}\ .$$
C'est le nombre minimum d'\'eclatements n\'ecessaires pour s\'eparer des branches index\'ees par $T$ (ou encore, $\alpha(T)-1$ est la hauteur de la 
bifurcation correspondant \`a $T$).\\

Pour $T \in{\cal B}$, nous notons \'egalement~:\\
$$\alpha'(T)=\left\{ \begin{array}{l}
 0 \mbox{ si }T = \{1,2,\ldots,n\}\\
\mbox{sup}\{ m_{i,k} \ ;\  i\in T ,\, k\notin T \} \;$sinon.$                                         
\end{array}\right. $$

Lorsque $\alpha'(T)$ est strictement positif, il existe un plus petit majorant strict $T'\supset T$
dans $\cal B$ d\'efini par~: 
$$ T' = \{ k \, ;\, \exists \, i \in T \:
\mbox{ tel que}\,:\ m_{i,k}\geq \alpha'(T) \}$$
et on a $\alpha (T') = \alpha'(T)$. 
La bifurcation correspondant \`a $T'$ est la bifurcation
qui pr\'ec\`ede celle correspondant \`a $T$ dans l'arbre d'\'eclatements, et 
$\gamma (T) = \alpha (T) - \alpha' (T)$ est la longueur de la branche 
entre ces deux bifurcations (ou encore la diff\'erence des hauteurs des
bifurcations $T$ et $T'$).\\

Lorsque $\alpha'(T) = 0$, $T$ est une partie maximale 
de $\cal B$, et elle  correspond \`a une premi\`ere bifurcation de l'arbre. 
Dans tous les cas, nous posons :
$$\gamma (T) = \alpha (T) - \alpha'(T)\ .$$
Constatons enfin que l'ensemble $\cal B$ des parties de 
$\{1,2,\ldots,n\}$ v\'erifie les propri\'et\'es suivantes :

\indent - si $T \in{\cal B}$ alors $T$ poss\`ede au moins deux \'el\'ements ;\\
\indent - si $T_1,T_2\in{\cal B}$ alors $T_1 \bigcap T_2$ est vide, 
 ou $T_1 \subset T_2,$ ou $T_2 \subset T_1$\\
et la fonction $\gamma$ d\'efinie sur $\cal B$ est \`a valeurs dans les entiers naturels.\\

R\'eciproquement, \`a partir de $(\cal B,\gamma )$, nous pouvons 
reconstituer la famille des multiplicit\'es
$\{m_{i,j}\ ;\ i\neq j\}$ en posant :
$$m_{i,j}=\left\{ \begin{array}{l}
 0 \mbox{ si $\{i,j\}$ n'est inclu dans aucun \'el\'ement de $\cal B$} \\
 \sum_{\{i,j\}\subset T ,\, T\in \cal B}\; \gamma (T) \mbox{ sinon.} 
\end{array}\right. $$

\subsubsection{Le cas g\'en\'eral}

Consid\'erons maintenant un bouquet \`a $n$ branches. Par analogie
avec le cas pr\'ec\'edent, nous dirons qu'une
partie $T$ de $T_0=\{1,2,\ldots,n\}$ est un {\em rameau} associ\'e au 
bouquet si :

\indent - $T$ poss\`ede au moins deux \'el\'ements ;\\
\indent - pour tout $(i,j,k) \in {T^2}\times  T_0$ avec $i\neq j$
et $k \notin T$ : $m_{i,j} \, > \, m_{i,k} = m_{j,k}$.

\stl

Avec cette d\'efinition, la partie totale $T_0$ est donc 
toujours un rameau. De plus, les notions de
rameau et de bifurcation sont les mêmes pour un bouquet 
de courbes lisses - mise \`a part \'eventuellement la partie totale
(voir la remarque \ref{BouR}). En particulier, la famille $\cal R$ 
des rameaux v\'erifie les propri\'et\'es suivantes :

\indent - un rameau $T\in{\cal R}$ poss\`ede au moins deux \'el\'ements ;\\
\indent - si $T_1,T_2\in{\cal R}$ alors $T_1 \bigcap T_2$ est vide, 
ou $T_1 \subset T_2$, ou $T_2 \subset T_1$.\\

Pour $T \in {\cal R}$, nous posons comme pr\'ec\'edemment :
$$\alpha(T) = \mbox{inf}\{ m_{i,j} \ ;\ i, j\in T, i\neq j \}.$$ 
Constatons que la fonction 
$\alpha :\, {\cal R}\longrightarrow {\bf Q}$ est strictement d\'ecroissante.
Lorsque $T \in {\cal R}$ n'est pas la partie totale, $T$ poss\`ede un plus petit majorant strict $T'$ dans $\cal R$, et on pose : 
$\gamma (T) = \alpha (T) - \alpha (T')$. Pour la partie totale, 
on pose : $\gamma (T_0) = \alpha (T_0) = 
\mbox{inf}\{ m_{i,j}\ ;\ i\neq j \}$. 
Nous avons alors de nouveau, pour des indices $i\neq j$ distincts :
\begin{equation} \label{multinter}
m_{i,j} = \sum_{\{i,j\}\subset T , T\in \cal R}\; \gamma (T) 
\end{equation}
 
Nous dirons que $({\cal R},\gamma )$ est {\it l'arbre associ\'e au bouquet.}

\begin{defi}
 Une famille ${\cal R}$ de parties de $T_0=\{1,2,\ldots,n\}$ est une 
{\em esp\`ece d'arbre \`a $n$ branches} si elle v\'erifie les conditions :

 - la partie totale $T_0$ appartient \`a ${\cal R}$ ;

 - toute partie $T \in{\cal R}$ a au moins deux \'el\'ements ;

 - si $T_1,T_2  \in{\cal R}$, alors  $T_1 \bigcap T_2 $ est vide, 
ou $ T_1 \subset T_2,$ ou $T_2 \subset T_1$.
\end{defi}

\begin{defi}
Un {\em arbre} est un couple $({\cal R},\gamma )$ form\'e 
d'une esp\`ece d'arbre ${\cal R}$ \`a $n$ branches et 
d'une fonction $\gamma$ sur ${\cal R}$ 
\`a valeurs rationnelles, avec 
$ \gamma (T) > 0 $ pour $ T \neq T_0$.
\end{defi}

Nous pouvons dire que $\gamma$ donne l'altitude du premier
rameau et les dimensions de l'arbre. On d\'emontre 
facilement qu'un arbre est l'arbre associ\'e \`a un bouquet, 
et ce bouquet est un bouquet de courbes m\'eromorphes 
(resp. lisses) lorsque $\gamma$ est \`a valeurs dans ${\bf Z}$ 
(resp. dans ${\bf N}$).

\begin{rema}{\em \label{BouR}
Dans le cas d'un bouquet de courbes lisses, l'arbre
$({\cal R},\gamma)$ contient la m\^eme information
que le couple $({\cal B},\gamma)$. En effet, ceux-ci
co\"{\i}ncident lorsque $T_0\in{\cal B}$. Et quand
$T_0\not\in{\cal B}$, au moins deux courbes du bouquet 
coupent l'axe des $y$ en des points distincts {\it i.e.}
$m_{i,j}=0$ pour un couple $(i,j)$ ; par suite, $\gamma(T_0)=0$.
}
\end{rema}

\begin{exem}{\em
  L'arbre associ\'e au bouquet de courbes lisses consid\'er\'e dans
l'exemple \ref{exrec} est de l'esp\`ece suivante :
$${\cal R} \; = \; \{\{1,2,3,4,5,6\} , \,\{2,3\}, \, \{4,5,6\} , \, 
\{5,6\}\}\ ,$$
et la fonction $\gamma$ associ\'ee est d\'efinie par :
$$\left\{ \begin{array}{l}
\gamma (\{1,2,3,4,5,6\}) \; = \; 0 \\
\gamma (\{2,3\}) \; = \; 1 \\
\gamma (\{4,5,6\}) \; = \; 1 \\
\gamma (\{5,6\}) \; = \; 2
\end{array}\right. $$}
\end{exem}

\subsection{L'alg\`ebre magique associ\'ee \`a une esp\`ece d'arbre}
\label{paralgmag}

\subsubsection{D\'ecomposition des matrices sym\'etriques de
${\cal M}ag^\star(n,{\bf Q})$}

\`A une partie $T$ de $T_0=\{1,2,\ldots,n\}$ de cardinal $n(T)$ 
au moins \'egal \`a deux, nous associons la matrice magique 
exceptionnelle $A(T)=(\alpha _{i,j})$ sym\'etrique, de somme nulle, 
 d\'efinie par :
$$\alpha_{i,j}=\left\{ \begin{array}{cl}
 -1 &\mbox{si }i\neq j\mbox{ et }\{i,j\}\subset T\ ; \\
 n(T)-1 &\mbox{si } i=j\in T\ ; \\
  0 &\mbox{sinon.}
\end{array}\right. $$

Cette notation est motiv\'ee par le fait suivant.

\begin{prop} \label{decompmatmag} 
Soit $A$ la matrice magique associ\'ee \`a un bouquet \`a 
$n$ branches, et 
$({\cal R},\gamma )$ l'arbre associ\'e \`a ce bouquet. Alors :
$$ A \;\; = \sum _{T\in {\cal R}} \gamma (T) A(T) \ .$$
\end{prop}

C'est une cons\'equence directe de la formule (\ref{multinter}) page
\pageref{multinter}.

\stl

Dans la suite, $I$  d\'esignera toujours la matrice identit\'e de taille $n$.

\begin{coro} Soit $A\, = \, (\alpha _{i,j})$ une matrice magique exceptionnelle, sym\'etrique,  \`a coefficients rationnels.
Il existe un arbre $({\cal R},\gamma )$ tel que :
$$ A \; = \; s(A)I \; +\; \sum _{T \in {\cal R}} \gamma (T) A(T) $$
avec $\gamma (T_0) =\mbox{\em{inf}} \{ - \alpha _{i,j} \,\,;\; i\neq j \}.$
\end{coro}

\begin{rema}{\em
 Ce corollaire reste valable lorsque la matrice est \`a coefficients 
dans $\bf Z$ ou $\bf R$, la fonction $\gamma $ \'etant alors 
\`a valeurs dans $\bf Z$ ou $\bf R$.}
\end{rema}

\subsubsection{L'alg\`ebre magique associ\'ee \`a une esp\`ece d'arbre}
\label{algmag}

\'Etant fix\'ee une esp\`ece d'arbre $\cal R$, nous \'etudions 
ici l'alg\`ebre engendr\'ee par les matrices $A(T)$, $T\in{\cal R}$. 
Donnons d'abord quelques propri\'et\'es des matrices $A(T)$.

\begin{form} {\em \label{formucomred}
Soit $T_1$ et $T_2$ deux parties de $T_0 = \{1,2,\ldots,n\}$ de cardinal
$n(T_1)$ et  $n(T_2)$ au moins \'egal \`a deux. Alors :

 - si $T_1 \bigcap T_2 =\emptyset \; : \; A(T_1)A(T_2) \, = \, A(T_2)A(T_1) \, = 0$ ;

 - si $T_1 \subset T_2 \; : \; A(T_1)A(T_2) \, = \, A(T_2)A(T_1) \, = n(T_2)A(T_1)$. }
\end{form}

Afin de diagonaliser $A(T)$, fixons maintenant quelques notations. 

\begin{nota}{\em \label{notared}
Soit $E = {\bf C}^n$ et $\{e_1, \ldots , e_n \}$ sa base canonique. 
Pour une partie $T$ de $T_0 = \{1,2,\ldots,n\}$ de cardinal $n(T)$, 
notons :
$E(T)={\bigoplus}_{i \in T} {\bf C} e_i$,   
$F(T)=\left\{ \sum _{i \in T} u_i e_i \; ; \; \sum _{i \in T } u_i = 0\right\}$ et $\omega (T) \, = \, \sum _{i \in T} e_i$. }
\end{nota}

Nous convenons d'identifier un endomorphisme de $E$ 
et sa matrice dans la base canonique. 
Constatons que : $A(T)(e_j)=-\omega(T)+n(T)e_j$ 
si $j\in T$, et $A(T)(e_j)=0$ sinon.
Il vient alors ais\'ement :
$$ \begin{array}{l}
E = E(T^c) \oplus  F(T)  \oplus {\bf C} \omega (T) \\
\ker A(T) =  E(T^c)  \oplus  {\bf C} \omega (T) \\
\ker\;(A(T) - n(T)I) = F(T)
\end{array}$$
o\`u $T^c = T_0 - T$ d\'esigne le compl\'ementaire de $T$. 
En particulier,  $A(T)$ est diagonalisable dans $\bf Q$ ; 
ses valeurs propres sont $0$ et $n(T)$, et les 
sous-espaces propres correspondants sont
respectivement $E(T^c)\oplus {\bf C}\,\omega (T)$ et $F(T)$.

\stl

Donnons deux cons\'equences de ces r\'esultats.

\begin{prop} Soit ${\cal R}$ une esp\`ece d'arbre.

(i)  Le $\bf Q$-espace vectoriel de base $\{ A(T) \; ; \; T 
\in {\cal R}\}$ est une $\bf Q$-alg\`ebre commutative de matrices magiques sym\'etriques, de somme nulle, diagonalisables.

(ii) Le $\bf Q$-espace vectoriel de base 
 $\{ I \} \bigcup \{A(T) \; ; \; T 
\in {\cal R}\}$ est une $\bf Q$-alg\`ebre commutative de matrices magiques sym\'etriques, diagonalisables.
\end{prop}

Nous dirons que la ${\bf Q}$-alg\`ebre d\'efinie au $(i)$ est la
{\em $\bf Q$-alg\`ebre magique associ\'ee} \`a l'esp\`ece d'arbre 
${\cal R}$. 

 \stl

\begin{demo}
Constatons que toutes les matrices $A(T)$, $T\in{\cal R}$, sont 
diagonalisables dans une m\^eme base (puisque elles sont diagonalisables
et commutent entre elles.)
  
Pour d\'emontrer que les matrices donn\'ees sont ind\'ependantes, 
nous v\'erifions d'abord que 
$e_1 + \cdots + e_n $ est un vecteur propre de $I$ qui 
est dans le noyau de toutes les matrices $A(T)$. 
Il reste alors \`a montrer
que pour tout \'el\'ement maximal $T_1$ d'une sous-famille 
${\cal R'}\subset{\cal R}$, il existe un vecteur propre de $A(T_1)$ 
associ\'e \`a la valeur propre $n(T_1)$ et qui soit dans le noyau 
de $A(T)$ pour $T \in {\cal R'}$, $T \neq T_1$. Cela 
r\'esulte ais\'ement du fait suivant : pour tout \'el\'ement
$S\in{\cal R}$, on peut construire un vecteur $v_S\in E$ qui soit
un vecteur propre de $A(T)$, $T\in {\cal R}$, associ\'e \`a 
la valeur propre $n(T)$ lorsque $T\supset S$ et \`a la valeur
propre $0$ sinon.

\stl

Lorsque $S \in {\cal R}$ est minimal, on constate en effet que tout vecteur 
non nul de $F(S)$ convient (et $F(S)\neq 0$
 puisque $n(S)$ est sup\'erieur ou \'egal \`a $2$, ${\cal R}$ \'etant
une esp\`ece d'arbre). Lorsque $S \in {\cal R}$ n'est pas minimal, notons $S_1 , \ldots , S_\ell$
les minorants stricts maximaux de $S$ dans ${\cal R}$ et $R$ le
compl\'ement dans $S$ : $S$ est la r\'eunion disjointe :
$S_1 \bigcup \cdots \bigcup S_\ell \bigcup R $. 
Si $R$ n'est pas vide\footnote{C'est en particulier le cas lorsque $\ell=1$}, 
nous avons la d\'ecomposition :
$$F(S) = F(R)\oplus(\bigoplus_{i=1}^\ell  F(S_i))  
\oplus (\bigoplus_{i=1}^\ell {\bf C}(\omega (S_i)
-\frac{n(S_i)}{n(S)}\omega (S)))\ .$$
Si au contraire $R$ est vide, nous avons :
$$F(S) = (\bigoplus_{i=1}^\ell  F(S_i)) 
\oplus (\bigoplus_{i=1}^{\ell-1} {\bf C}(\omega (S_i)
-\frac{n(S_i)}{n(S)}\omega (S)))\ .$$
Dans les deux cas, le vecteur $v_S=\omega (S_1)-(n(S_1)/n(S))\omega (S)$ 
convient.
\end{demo}

\begin{prop}
   Soit ${\cal R}$ une esp\`ece d'arbre et 
$A= \sum _{T\in {\cal R}} \lambda (T) A(T)$
une matrice de la ${\bf Q}$-alg\`ebre magique associ\'ee. 
Alors $A$ est diagonalisable et ses valeurs propres
sont $0$ et :
$$ \; \sum _{T\in {\cal R},\;T \supset T_1} 
\lambda(T) n(T) \; \; , \; \; T_1 \in {\cal R} \ .$$
\end{prop}

Explicitons le cas particulier de la matrice magique 
associ\'ee \`a un bouquet passant par 
l'origine (D\'efinition \ref{bouqpass}).
Dans ce cas, $T_0 = \{1,2,\ldots,n\} \;$ est 
la bifurcation maximale de l'arbre associ\'e au bouquet ; en particulier,
$\gamma(T_0)>0$ et le noyau de la matrice magique correspondante est 
${\bf C}\omega$ avec
 $\omega = \omega (T_0) = e_1 + \cdots + e_n $. Nous avons alors :
 $E=F\oplus{\bf C}\,\omega$ o\`u $F=F(T_0)
=\{\sum _{i=1}^n u_i e_i \; ; \; \sum _{i=1}^n u_i=0\}$.

\begin{coro} \label{AFcourbliss}
Soit $A$ la matrice magique associ\'ee \`a un bouquet de branches 
(resp. de courbes lisses) 
passant par l'origine et $({\cal R},\gamma )$
l'arbre associ\'e. Alors $A$ induit 
un automorphisme de $F$ dont les valeurs propres sont les 
rationnels (resp. les entiers)
strictement positifs~:
$$\sum _{T\in {\cal R},\;T \supset
 T_1} \gamma (T) n(T) \; \; , \; \; T_1 \in {\cal R}\ .$$
\end{coro}

\begin{exem}{\em \label{remAFcourbliss}
 Dans le cas du bouquet consid\'er\'e \`a l'exemple \ref{exrec},
les valeurs propres de la matrice magique associ\'ee sont : 
$0, \, 2, \, 3, 7$.}
\end{exem}

\subsection{La matrice magique g\'en\'eralis\'ee associ\'ee 
\`a un bouquet de branches avec multiplicit\'es}
\label{magikgen}

\`A un bouquet de branches avec multiplicit\'es  
(D\'efinition \ref{defbouqmult}), nous associons
 ici une matrice carr\'ee qui g\'en\'eralise la
notion de matrice magique d'un bouquet de branches.

\begin{defi} \label{matmagen}
   La {\em matrice magique g\'en\'eralis\'ee} associ\'ee au bouquet de
branches avec multiplicit\'es d\'efini par le polyn\^ome 
 $\prod_{i=1}^n  (y - a_i)^{\mu_i}\in\overline{K}[y]$
est la matrice carr\'ee \`a $n$ lignes et $n$ colonnes 
$A_\mu=(\alpha_{i,j})$, \`a coefficients rationnels d\'efinis par :
$\alpha_{i,j} = - \mu_i m_{i,j} $ pour $i\not=j$, et
 $\alpha_{i,i} =  \sum_{k \not=i} \, \mu_k m_{k,i}$.
\end{defi}

En g\'en\'eral, cette matrice n'est bien s\^ur ni magique,
ni sym\'etrique. Toutefois, la somme des coefficients de
chacune de ses colonnes est nulle. Nous avons aussi
l'identit\'e matricielle~:

$$A_\mu=\left( \begin{array}{ccc}
  \mu_1 &\   & 0\\
 & \ddots & \\
0 &\  & \mu_n
\end{array} \right)
\cdot A+
\left(\begin{array}{ccc}
\sum_j(\mu_j-\mu_1)m_{1,j} &\  & 0  \\
 \ & \ddots & \ \\ 
0 & \ & \sum_j(\mu_j-\mu_n)m_{n,j}
\end{array}
\right) $$
o\`u $A$ est la matrice magique associ\'ee au bouquet de $n$
branches sous-jacent. Cette relation  permet
d'\'etendre \`a la matrice $A_\mu$ les r\'esultats sur $A$ 
obtenus au paragraphe \ref{paralgmag}.

\stl

Soit $\mu = (\mu_1,\ldots ,\mu_n)$ une famille 
d'entiers naturels non nuls. \`A une partie $T$ de $\{1,2,\ldots,n\}$ 
de cardinal $n(T)$ au moins \'egal \`a deux, nous
associons la matrice $\; A_\mu (T)=(\alpha_{i,j})$ 
d\'efinie par~:
$$\alpha_{i,j}=\left\{ \begin{array}{cl}
   -\mu_i & \mbox{si } i\neq j \mbox{ et } \{i,j\}\subset T\ ; \\
\sum_{k\neq i , k\in T}\, \mu_k & \mbox{si } i=j\in T\ ; \\
 0 &\mbox{sinon.}
\end{array}\right. $$

\begin{prop} 
Soit $A_\mu$ la matrice magique g\'en\'eralis\'ee  associ\'ee 
\`a un bouquet de branches avec multiplicit\'es, et 
$({\cal R},\gamma )$ l'arbre associ\'e au  bouquet de
branches sous-jacent. Alors :
$$ A_\mu \;\; = \sum _{T\in {\cal R}} \gamma (T) A_\mu(T) \ .$$
\end{prop}

C'est une cons\'equence directe de la proposition  
\ref{decompmatmag} et de l'identit\'e matricielle reliant
$A_\mu$ et $A$. 
Constatons que le formulaire
\ref{formucomred}  se g\'en\'eralise aussi.

\begin{form}{\em
Soit $T_1$ et $T_2$ deux parties de $T_0=\{1,2,\ldots,n\}$ de
cardinal $n(T_1)$ et  $n(T_2)$ au moins \'egal \`a $2$. Alors :

 - si $T_1 \bigcap T_2 = \emptyset$ :  $A_\mu (T_1)A_\mu (T_2) = A_\mu (T_2)A_\mu (T_1) = 0$ ;

 - si $T_1 \subset T_2$ : $A_\mu (T_1)A_\mu (T_2) = A_\mu (T_2)A_\mu (T_1)
  =  \sigma_{\mu}(T_2)A_\mu (T_1)$ 

\noindent o\`u pour toute partie $T$ de $T_0$, 
$\sigma_\mu(T)$ d\'esigne la somme $\sum_{i\in T}\mu_i$.}
\end{form}

De plus, avec les notations \ref{notared}, 
 le noyau de $A_\mu (T)$ est $E(T^c)\bigoplus {\bf C}\,
\omega_\mu(T)$ avec $\omega_\mu(T)=\sum_{i\in T}\mu_i e_i$,
et le sous-espace propre associ\'e \`a la valeur propre 
$\sigma_{\mu}(T)$ est $F(T)$. En particulier, 
$A_\mu (T)$ est diagonalisable.

\begin{prop} \label{ValpropAmu}
Soit $A_\mu$ la matrice magique g\'en\'eralis\'ee
associ\'ee \`a un bouquet de branches avec
multiplicit\'es passant par l'origine et $({\cal R},\gamma )$
l'arbre associ\'e. Alors $A_\mu$ induit 
un automorphisme de $F=F(T_0)$ 
dont les valeurs propres sont les rationnels 
strictement positifs~:
$$\sum _{T\in {\cal R},\;T \supset
 T_1} \gamma (T)  \sigma_{\mu}(T) \; \; , \; \; T_1 \in {\cal R}\ .$$
\end{prop}

Ce r\'esultat s'obtient de la même mani\`ere que dans le cas r\'eduit ;
 nous ne recopions pas, ici, sa d\'emonstration d\'etaill\'ee.

\subsection{Le cas particulier des germes irr\'eductibles}

Nous d\'etaillons ici les constructions et
r\'esultats pr\'ec\'edents dans le cas d'un bouquet de branches
associ\'e \`a une courbe irr\'eductible, transverse \`a l'axe des
$y$, et d\'efinie par un polyn\^ome (irr\'eductible) 
$f\in {\bf C}\{x\}[y]$ unitaire de degr\'e $n\geq 2$ en $y$.

\stl

D'apr\`es le th\'eor\`eme de Newton-Puiseux, $f$ se factorise
dans ${\bf C}[[x^{1/n}]][y]$ :
$$ f=\prod_{i=1}^n (y-\phi(\xi^ix^{1/n}))$$ 
o\`u  $\phi(t)=\sum_{i\geq n} u_i t^i\in{\bf C}[[t]]$ et 
$\xi\in{\bf C}$ est une racine primitive $n$-\`eme de 
l'unit\'e ;
 nous posons alors $a_i=\phi(\xi^i x^{1/n})$ pour $1\leq i\leq n$. 
Les multiplicit\'es d'intersection $m_{i,j}$ s'expriment alors
simplement en fonction des exposants caract\'eristiques 
$(\beta_0,\ldots, \beta_g)$ de la courbe d\'efinie par $f$.
Rappelons qu'ils sont d\'efinis \`a partir de $\phi$ de la
fa\c con suivante : $\beta_0=n$ et pour $k\geq 1$, 
$\beta_k=\mbox{inf}\{i\, ; \, u_i\not=0 
\mbox{ et p.g.c.d}(i,\beta_0,\ldots,\beta_{k-1})<\mbox{p.g.c.d}
(\beta_0,\ldots, \beta_{k-1})\}$. On pose alors : 
$\epsilon_k=\mbox{p.g.c.d}\{\beta_0,\ldots,\beta_k\}$, $0\leq k\leq g$, 
et $n_k=\epsilon_{k-1}/\epsilon_k$, $1\leq k\leq g$. En particulier,
nous avons $\epsilon_0=n$, $\epsilon_k=n_{k+1}\cdots n_g$
pour $0\leq k<g$, $\epsilon_g=1$, $n=n_1\cdots n_g$,
et $\phi(t)$ se d\'ecompose en une somme :
$$\phi(t)=t^n\phi_0(t^n)+t^{\beta_1}\phi_1(t^{\epsilon_1})+\cdots+
t^{\beta_g}\phi_g(t^{\epsilon_g})$$
o\`u $\phi_k(t)\in{\bf C}[[t]]$ est inversible pour $1\leq k\leq g$.
Les multiplicit\'es d'intersection $m_{i,j}$, $i\not=j$, sont alors
donn\'ees par :
$$m_{i,j}=\left\{
  \begin{array}{rcl}
    \beta_1/n &\mbox{lorsque} & |i-j| \mbox{ n'est pas un multiple de $n_1$\ ;}\\
    \beta_k/n &\mbox{lorsque} & |i-j| \mbox{ est un multiple de 
       $n_1\cdots n_{k-1}$} \\
      & & \mbox{ qui n'est pas divisible par $n_1\cdots n_k$.}
   \end{array} \right.$$  

Remarquons alors que la matrice magique $A=(\alpha_{i,j})$ 
associ\'ee au bouquet v\'erifie la propri\'et\'e suivante :
les coefficients $\alpha_{i,j}$ sont constants le long des
`petites diagonales' $i=j\pm\ell$, $1\leq\ell\leq n-1$. De plus,
un calcul facile montre que les termes diagonaux $\alpha_{i,i}$ sont
tous \'egaux \`a $\sum_{k=1}^g(\epsilon_{k-1}-\epsilon_k)(\beta_k/n)$.

\begin{exem}{\em
Si $n=6$ et $\phi(t)=t^8+t^9$, alors $\beta_0=6$, $\beta_1=8$,
$\beta_2=9$, $\epsilon_0=6$, $\epsilon_1=2$ et $\epsilon_2=1$.
La matrice magique associ\'ee au bouquet est alors :
$$\left(\begin{array}{cccccc}
    41/6 & -4/3 & -4/3 & -3/2 & -4/3 & -4/3 \\
    -4/3 & 41/6 & -4/3 & -4/3 & -3/2 & -4/3 \\
    -4/3 & -4/3 & 41/6 & -4/3 & -4/3 & -3/2  \\
    -3/2 & -4/3 & -4/3 & 41/6 & -4/3 & -4/3  \\
    -4/3 & -3/2 & -4/3 & -4/3 & 41/6 & -4/3   \\
    -4/3 & -4/3 & -3/2 & -4/3 & -4/3 & 41/6 
    \end{array}\right)$$}
\end{exem}

Nous allons maintenant expliciter l'arbre associ\'e au bouquet.
Lorsque $g=1$, l'esp\`ece d'arbre ${\cal R}$
est r\'eduit \`a la partie totale $T_0$ et $\gamma(T_0)=\beta_1/n$.
Supposons maintenant que $g$ soit au moins \'egal \`a $2$. La
famille ${\cal R}$ est alors constitu\'ee de $T_0$ et des parties
$T_{i_1,\ldots,i_r}$, $1\leq r\leq g-1$, $1\leq i_1\leq n_1$ et
$0\leq i_k \leq n_k-1$ pour $2\leq k\leq r$, d\'efinies par :
$$T_{i_1,\ldots,i_r}=\left\{i_1+\sum_{k=2}^r i_k\times n_1\cdots
n_{k-1} +\ell\times n_1 \cdots  n_r\ ; \ 0\leq \ell 
\leq \epsilon_{r}-1  \right\}\ .$$
En d'autres termes, $T_{i_1,\ldots,i_r}$ est l'ensemble des indices
$i\in T_0$ tels que pour tout $1\leq j\leq r$, le reste de la division
euclidienne de $i$ par $n_1\cdots n_j$ soit \'egal \`a
$i_1+\sum_{k=2}^j i_k\times n_1\cdots n_{j-1}$. Ainsi 
$n(T_{i_1,\ldots,i_r})=\epsilon_r$, $T_0$ est la r\'eunion disjointe
des parties $T_1$,\ldots, $T_{n_1}$, et pour $g\geq 3$, $1\leq r\leq g-2$,
$T_{i_1,\ldots,i_r}$ est la r\'eunion disjointe de 
$T_{i_1,\ldots,i_r,0},\ldots, T_{i_1,\ldots, i_r,n_{r+1}-1}$.

La fonction $\gamma$ est d\'efinie par $\gamma(T_0)=\beta_1/n$ et
$\gamma(T_{i_1,\ldots,i_r})=(\beta_{r+1}-\beta_r)/n$ pour $1\leq r\leq g-1$
(ind\'ependamment de la valeur des indices $i_1,\ldots, i_r$).

\stl

Il r\'esulte alors directement du corollaire \ref{AFcourbliss} 
que l'endomorphisme de $F$ induit par la matrice magique $A$ a 
pour valeurs propres : 
$$\beta_1,\ \beta_1+(\beta_2-\beta_1)\epsilon_1/n,\ \ldots,\  
\beta_1 +(\beta_2-\beta_1)\epsilon_1/n +\cdots +(\beta_g-\beta_{g-1})\epsilon_{g-1}/n\ .$$

\section{Division par l'id\'eal jacobien}

\subsection{Pr\'eliminaires}
\label{parcalcwf}
Soit $f=\prod_{i=1}^{n} (y - a_{i})\in\overline{K}[y]$  un bouquet \`a $n$ branches passant par l'origine 
(D\'efinition \ref{bouqpass}) ; en particulier, les 
polyn\^omes $f'_x$ et $f'_y$ sont premiers 
entre eux dans $\overline{K}[y]$.  

\stl

 Notons ${\cal E} = K[y]^{(n)}$ et $\overline{\cal E} = \overline{K}[y]^{(n)}$ 
les espaces de polynômes de degr\'e en $y$ strictement
inf\'erieur \`a $n$, \`a coefficients dans $K$ et $\overline{K}$ 
respectivement, 
 ${\cal F} = K[y]^{(n-1)}$ et $\overline{\cal F} = \overline{K}[y]^{(n-1)}$ 
les sous-espaces de polynômes de degr\'e strictement inf\'erieur \`a $n-1$.
Remarquons alors que la famille de polyn\^omes~:
$$\varepsilon_i = \prod_{j \neq i} (y - a_{j}) \ ,\; \ 1\leq i \leq n\, , $$
forme une base du $\overline{K}$-espace vectoriel $\overline{\cal E}$. 
En particulier,  l'espace  $\overline{\cal F}$ s'identifie \`a l'hyperplan
de $\overline{\cal E}$~: 
$\{\sum_{i=1}^{n}u_i \varepsilon_i \; ;\; u_i\in{\overline K} \; , \; \sum_{i=1}^{n} u_i = 0 \, \}$, et nous avons les d\'ecompositions :  
${\cal E} = {\cal F} \oplus {K}\omega$, 
 $\overline{\cal E} = \overline{\cal F} \oplus \overline{K}\omega$ avec 
$\omega = \sum_{i=1}^{n} \varepsilon_i = f'_y$.

\stl

Comme $f'_x$ et $f'_y$ sont premiers entre eux, l'\'equation :
$$u f'_x +  v f'_y = w $$
admet une unique solution $(u,v)$ dans 
$\overline{\cal F} \times \overline{K}[y]$ pour tout 
$w \in \overline{K}[y]$. Lorsque
$w=\tilde{w}f$ avec $\tilde{w} \in \overline{\cal F}$, cette 
solution est dans $\overline{\cal F} \times \overline{\cal E}$ 
et elle s'obtient par r\'esolution d'un syst\`eme de Cramer. Dans
les paragraphes qui suivent, nous allons expliciter cette solution en 
travaillant dans la base $(\varepsilon_1, \ldots , \varepsilon_n)$.
\`A cette fin,  introduisons maintenant une filtration naturelle sur 
les espaces  $\overline{\cal E}$, $\overline{\cal F}$, $\cal E$ et 
$\cal F$.

\stl

 Consid\'erons  $(\overline{K}_r)_{r\in{\bf Q}}$ la filtration d\'ecroissante 
du corps $\overline{K}$ d\'efinie par : 
$\overline{K}_r=\{a\in\overline{K}\, ;\, \nu(a)\geq r \}$ pour
$r\in {\bf Q}$. Nous munissons alors les espaces $\overline{\cal E}$ et 
$\overline{\cal F}$ des filtrations :
$$\overline{\cal E}_r=\bigoplus_{i=1}^n\overline{K}_r\varepsilon_i\ ,\ \; 
  \overline{\cal F}_r=\overline{\cal E}_r \cap \overline{\cal F}
=\bigoplus_{i=1}^{n-1}\overline{K}_r(\varepsilon_i-\varepsilon_n)$$
d\'efinies par le `poids' : 
$$\mbox{val}( \sum_{i=1}^{n} u_i \varepsilon_i )= \mbox{inf} \{ \nu(u_1), \ldots , \nu(u_n) \}\ .$$
Ces filtrations induisent bien s\^ur \`a leur tour des filtrations sur
${\cal E}$ et ${\cal F}$ en posant, pour tout $r\in{\bf Q}$ : 
${\cal E}_r=\overline{\cal E}_r\cap {\cal E}$ et
${\cal F}_r=\overline{\cal F}_r\cap {\cal F}$. Pour 
$r\in{\bf Q}$, notons encore $\overline{\cal E}_{>r}$ 
(resp. $\overline{\cal F}_{>r}$, 
${\cal E}_{>r}$, ${\cal F}_{>r}$) le sous-espace de $\overline{\cal E}_{r}$ 
(resp. $\overline{\cal F}_r$, ${\cal E}_r$, ${\cal F}_r$) 
form\'e des \'el\'ements de valuation strictement sup\'erieure \`a $r$.

\begin{rema} {\em \label{Ezero}
(i) Les polyn\^omes $\varepsilon_i$, $1\leq i\leq n$, 
 interviennent bien s\^ur dans la formule d'interpolation de Lagrange ;
en particulier, si $u\in \overline{\cal E}$ alors 
$\sum_{i=1}^n(u(x,a_i)/\varepsilon_i(a_i))\varepsilon_i$ est 
sa d\'ecomposition dans la base $(\varepsilon_1,\ldots,\varepsilon_n)$.

(ii) Il est ais\'e de constater  que  ${\cal E}_0$ est 
inclu dans ${\bf C}[[x]][y]$
lorsque le bouquet d\'efini par $f$ passe par l'origine.

(iii) Avec les notations introduites au paragraphe 
\ref{algmag}, le terme de degr\'e $r$ du gradu\'e associ\'e \`a la filtration 
$(\overline{\cal E}_r)_{r \in {\bf Q}}$ de $\overline{\cal E}$, 
$ \mbox{gr}_r\overline{\cal E}$,  s'identifie \`a $x^r E$ 
pour $r\in{\bf Q}$ ;  de m\^eme, $ \mbox{gr}_r\overline{\cal F}$ 
s'identifie \`a $x^r F$.}
\end{rema}

Pour un \'el\'ement $u$ de $\overline{\cal E}_r$, nous noterons enfin 
$\mbox{in}_r\,u\in\mbox{gr}_r\overline{\cal E}$ sa partie initiale 
en degr\'e $r$.

\newpage
\subsection{La r\'esolution de l'\'equation $ u f'_x  +  v f'_y  =  w f$}

Ce paragraphe est consacr\'e \`a l'\'etude de l'\'equation~:
\begin{equation} \label{eqwf}
u f'_x +  v f'_y = wf 
\end{equation}
lorsque $f$ d\'efinit un bouquet passant par l'origine,
avec ou sans multiplicit\'es.

\subsubsection{Le cas d'un bouquet \`a $n$ branches}

La r\'esolution de l'\'equation (\ref{eqwf}) s'appuie sur le r\'esultat
suivant.

\begin{lemm} \label{lemmprepawf}
  Soit $f\in\prod^{n}_{i=1}(y-a_i)\in\overline{K}[y]$ un polyn\^ome
r\'eduit de degr\'e $n\geq 2$. Soit $w=\sum_{i=1}^nw_i\varepsilon_i$
un \'el\'ement de $\overline{\cal F}$. Toute  solution 
$(u=\sum_{i=1}^n u_i\varepsilon_i,v)$
dans $\overline{\cal E}\times\overline{\cal E}$ de l'\'equation~:
$uf'_x+vf'_y=wf$ v\'erifie~: $v=\sum_{i=1}^n u_ia'_i\varepsilon_i$, et~:
$$ w_i = ( \sum_{j\neq i} \frac{a'_i - a'_j}{a_i - a_j})  u_i  - \sum_{j\neq i} (\frac{a'_i - a'_j}{a_i - a_j})  u_j \ ,\ \; 1\leq i\leq n\ .$$
\end{lemm}

\begin{demo}
Posons $v=\sum_{i=1}^{n} v_i \varepsilon_i$. Apr\`es division par $f^2$, 
l'\'equation (\ref{eqwf}) s'\'ecrit :
$$\sum_{i,j}\; \frac{-u_i a'_j + v_i}{(y-a_i)(y-a_j)} \; \; = \; \; \sum_{i=1}^n \; \frac{w_i}{y-a_i}\ .$$
Nous en d\'eduisons d'abord~: $v_i = u_i a'_i \ ,\ \; 1\leq i\leq n$. 
Par ailleurs, nous avons la d\'ecomposition :
\begin{equation} \label{fracrat}
 \frac{1}{(y-a_i)(y-a_j)}  =  \frac{1}{a_i - a_j} 
\left[ \frac{1}{y-a_i}  -  \frac{1}{y-a_j}\right] \ ,   
\end{equation}
pour tout $i\not=j$. En substituant dans le syst\`eme, nous obtenons alors :
$$\sum_{i\neq j} \frac{u_i (a'_i - a'_j)}{a_i - a_j}  
\left[ \frac{1}{y-a_i} - \frac{1}{y-a_j}\right] =  \sum_{i=1}^n \; \frac{w_i}{y-a_i}\ .$$
D'o\`u l'assertion.
\end{demo}
  
\begin{defi}
  Soit $f=\prod^n_{i=1}(y-a_i)\in\overline{K}[y]$ un polyn\^ome
r\'eduit de degr\'e $n\geq 2$. 
  La {\em matrice magique compl\`ete associ\'ee au bouquet} 
d\'efini par $f$ est la matrice magique  ${\cal A}$ de taille $n\times n$,  sym\'etrique, de somme nulle, \`a coefficients dans $\overline K$, de terme 
g\'en\'eral : $- (a'_i-a'_j)/(a_i-a_j)$ pour $1\leq i\not=j\leq n$.
\end{defi}

Si $\cal U$ et $\cal W$ sont les matrices colonnes form\'ees
respectivement par les coordonn\'ees de $u$ et $w$ dans
la base $(\varepsilon_1,\ldots,\varepsilon_n)$ de $\overline{\cal E}$,
nous avons alors l'identit\'e matricielle~: 
$${\cal W} = {\cal A}\cdot{\cal U}\ .$$
En particulier, on peut expliciter la solution de l'\'equation
(\ref{eqwf}) lorsque ce syst\`eme lin\'eaire est inversible.

\begin{prop} \label{ArestFauto}
   Soit $f=\prod^n_{i=1}(y-a_i)\in\overline{K}[y]$ un polyn\^ome 
r\'eduit d\'efinissant un bouquet passant par l'origine. La
matrice magique compl\`ete $\cal A$ associ\'ee induit un 
automorphisme de $\overline{\cal F}$ not\'e 
${\cal A}_{|\overline{\cal F}}$. De plus, le gradu\'e du morphisme
inverse $({\cal A}_{|\overline{\cal F}})^{-1}$ est de degr\'e $+1$ 
et s'identifie \`a $x(A_{|F})^{-1}$ o\`u $A_{|F}$ est la restriction 
\`a $F$ de l'endomorphisme de $E$ induit par la matrice 
magique $A$ associ\'ee au bouquet. 
\end{prop}

\begin{demo}
Pour $i \neq j $, nous posons : $ a_i - a_j \, = \, x^{m_{i,j}} \, c_{i,j}$, 
o\`u $c_{i,j}\in\overline{K}$ est de valuation nulle. 
Nous avons alors :
$$- \frac{a'_i - a'_j}{a_i - a_j}=- \frac{m_{i,j}}{x}-\frac{c'_{i,j}}{c_{i,j}}$$
Par suite, la matrice $\cal A$ se d\'ecompose en la somme :
$ {\cal A}= (1/x)A + {\cal A'}$,
o\`u ${\cal A'}$ est une matrice magique, sym\'etrique, 
de somme nulle, \`a coefficients dans $\overline{K}$,
 de valuation\footnote{Cette valuation est positive ou nulle 
lorsque nous avons affaire \`a un bouquet de branches 
m\'eromorphes.} strictement sup\'erieure \`a $-1$.

\stl

Comme la restriction de $A$ \`a $F$ est un endomorphisme inversible 
(Corollaire \ref{AFcourbliss}), la matrice  ${\cal A}$ induit donc
un automorphisme de $\overline{\cal F}$ d'inverse : 
$${\cal (A_{| \overline F})}^{-1} =  \left(  I  
+  x(A_{|\overline {\cal F}})^{-1} {\cal A'_{| \overline F}}\right)^{-1}
 x(A_{| F})^{-1} =  x(A_{|F})^{-1} + x  {\cal S}\ ,$$
o\`u ${\cal S}$ d\'esigne un endomorphisme de ${\overline {\cal F}}$ 
dont la matrice dans la base 
$(\varepsilon_1 - \varepsilon_n , \ldots , \varepsilon_{n-1} - \varepsilon_n)$ 
a ses coefficients de valuation strictement positive. 
Sous les identifications faites \`a la remarque \ref{Ezero}, 
le gradu\'e associ\'e au morphisme 
${\cal (A_{| \overline F})}^{-1}$ est bien de degr\'e $+1$, et :
$\mbox{gr}\,{\cal (A_{| \overline F})}^{-1} =  x  (A_{|F})^{-1}$.
\end{demo}

\begin{coro} \label{laresolwf}
  Soit $f=\prod^n_{i=1}(y-a_i)\in\overline{K}[y] $ un 
polyn\^ome r\'eduit 
d\'efinissant un bouquet passant par l'origine. Notons $A$ et 
${\cal A}$ les matrices magiques associ\'ees. Pour 
$w$ appartenant \`a $\overline{\cal F}_r$, 
la solution dans $\overline{\cal F} \times \overline{\cal E}$ 
de l'\'equation $ u f'_x +  v f'_y = w f $ est donn\'ee par :
$$\left\{  \begin{array}{ccccl}
    u  & = &  {\cal (A_{| \overline F})}^{-1}  w  &\in&  \overline{\cal F}_{r+1} \\ 
    v  & = & \sum_{i=1}^n u_i a'_i \varepsilon_i &\in& \overline{\cal E}_{>r}
\end{array}\right. $$
o\`u $u=\sum_{i=1}^n u_i \varepsilon_i$. 
En particulier~: $$\overline{\cal F}_{r}f \subset  \overline{\cal F}_{r+1} f'_x  +  
\overline{\cal E}_{>r} f'_y\ .$$
Plus pr\'ecis\'ement, si $\mbox{\em{in}}_r\,w= x^r W \in x^rE$, alors 
$\mbox{\em{in}}_{r+1}\,u = x^{r+1} (A_{|F})^{-1}\,W$. 
\end{coro}

\begin{rema}{\em \label{remareswf}
(i) Lorsque $w$ est de degr\'e sup\'erieur ou \'egal \`a $n-1$, la solution
$(u,v)$ dans $\overline{\cal F}\times \overline{K}[y]$ de l'\'equation
(\ref{eqwf}) se d\'eduit bien s\^ur du r\'esultat pr\'ec\'edent, en 
rempla\c cant $w$ par son reste dans la division par $f'_y$. En particulier, 
pour r\'esoudre l'\'equation $u f'_x + v f'_y=f^2$, on appliquera la proposition
pr\'ec\'edente  \`a $w=f-(y/n - \sigma/n^2)f'_y =  \sum_{i=1}^n (\sigma/n^2 - a_i/n)\varepsilon_i$ avec $\sigma = \sum_{i=1}^n a_i$.

(ii) Lorsque tous les \'el\'ements $a_i$ appartiennent au corps interm\'ediaire 
$K'={\bf C}[[x^{1/d}]][1/x]$ pour un certain entier naturel $d$, 
nous pouvons remplacer dans tous les calculs pr\'ec\'edents $\overline{K}$
par $K'$ ; les valuations sont alors \`a valeurs dans $(1/d){\bf Z}$. 
En particulier si $f$ est un bouquet de courbes lisses passant par l'origine, 
tout se passe dans le corps $K = {\bf C}((x))$.}
\end{rema}

Pr\'ecisons le corollaire pr\'ec\'edent dans le cas d'un polyn\^ome 
distingu\'e r\'eduit (Exemple \ref{expoldist}). 

\begin{coro} \label{cororesolwf}
 Soit $f\in {\bf C}[[x]][y]$ un polynôme distingu\'e, r\'eduit, de degr\'e 
$n\geq 2$. Notons $d\in {\bf N}$ le p.p.c.m des degr\'es de ses facteurs 
irr\'eductibles. 
Alors, pour tout $r \in (1/d){\bf Z}$ :
$${\cal F}_{r}f \subset  {\cal F}_{r+1} f'_x  +  {\cal E}_{>r} f'_y\ .$$
\end{coro}

\begin{demo}	
En effet, pour tout \'el\'ement $w \in {\cal F}$, la solution 
$(u,v)$ du syst\`eme $K$-lin\'eaire $u f'_x +  v f'_y  = w f$ 
est bien dans ${\cal F} \times {\cal E}$.
\end{demo}

Comme exemple d'application, prenons  $w=f-( y/n- \sigma/n^2)f'_y = 
\sum_{i=1}^n (\sigma/{n^2} - a_i/n)\varepsilon_i$, qui appartient 
 \`a ${\cal F}_{>0}$. Nous obtenons alors :
$$ wf \in {\cal F}_{>1} f'_x  + {\cal E}_{>0} f'_y $$
En se servant de la remarque \ref{Ezero}, nous en d\'eduisons que 
$f^2$ appartient \`a l'id\'eal jacobien. 
C'est l\`a l'id\'ee de la premi\`ere d\'emonstration de cette 
appartenance dans un manuscrit 
non publi\'e du premier auteur
(\cite{B}). Pour g\'en\'eraliser ce r\'esultat \`a un germe 
de courbe \`a singularit\'e isol\'ee, c'est-\`a-dire lorsque 
$f$ est le produit d'un polynôme distingu\'e r\'eduit par une
unit\'e, il convient d'it\'erer le corollaire et d'utiliser
le th\'eor\`eme d'annulation de Krull.

\subsubsection{Le cas d'un bouquet avec multiplicit\'es}

\label{resolwfnonred}

Dans ce paragraphe, nous nous int\'eressons \`a la
r\'esolution de l'\'equation~:
\begin{eqnarray} \label{eqwfmult}
u (f_{\mu}) '_x + v (f_{\mu})'_y = w f_{\mu} 
\end{eqnarray}
lorsque le polyn\^ome $ f_\mu = \prod_{i=1}^n (y - a_i)^{\mu_i}\in
\overline{K}[y]$ d\'efinit un bouquet de $n$ branches
avec multiplicit\'es, passant par l'origine 
(D\'efinition \ref{defbouqmult}).  Comme 
$(f_{\mu})'_x = g \prod_{i=1}^n (y - a_i)^{\mu_i - 1}$ et 
$(f_{\mu})'_y = h \prod_{i=1}^n (y - a_i)^{\mu_i - 1}$ avec
 $g$ et $h$ 
premiers entre eux dans ${\overline K}$[$y$], il y a encore 
existence et unicit\'e d'un couple solution 
$(u,v) \in \overline{\cal F} \times \overline{\cal E}$ 
pour l'\'equation (\ref{eqwfmult}) lorsque $w \in \overline{\cal F}$. 
La r\'esolution de cette \'equation fait aussi intervenir
une matrice carr\'ee \`a coefficients dans $\overline{K}$.

\begin{defi}
  Soit $f_\mu=\prod^n_{i=1}(y-a_i)^{\mu_i}\in\overline{K}[y]$ 
un polyn\^ome d\'efinissant un bouquet de $n$-branches
avec multiplicit\'es. Sa {\em matrice magique compl\`ete associ\'ee} 
est la matrice  ${\cal A}_\mu$ de taille $n\times n$  dont la
somme des termes de chacune de ses colonnes est nulle, 
\`a coefficients dans $\overline K$, de terme 
g\'en\'eral : $- \mu_i(a'_i-a'_j)/(a_i-a_j)$ pour $1\leq i\not=j\leq n$.
\end{defi}

\'Enon\c cons le pendant du corollaire \ref{laresolwf}.

\begin{prop}  \label{resolwfmultiple}
  Soit $f_{\mu}=\prod^n_{i=1}(y-a_i)^{\mu_i}\in\overline{K}[y]$ 
un polyn\^ome d\'efinissant un bouquet avec multiplicit\'es, 
passant par l'origine. Soit $A_\mu$ et 
${\cal A}_\mu$ les matrices  associ\'ees. Pour tout 
$w$ appartenant \`a $\overline{\cal F}_r$, 
la solution dans $\overline{\cal F} \times \overline{\cal E}$ 
de l'\'equation $ u (f_{\mu})'_x +  v (f_{\mu})'_y = w f_{\mu}$ 
est donn\'ee par :
$$\left\{  \begin{array}{ccccl}
    u  & = &  {\cal (A_{\mu| \overline F})}^{-1}  w  &\in&  \overline{\cal F}_{r+1} \\ 
    v  & = & \sum_{i=1}^n u_i a'_i \varepsilon_i &\in& \overline{\cal E}_{>r}
\end{array}\right. $$
o\`u $u=\sum_{i=1}^n u_i \varepsilon_i$. 
En particulier~: 
$$\overline{\cal F}_{r}f_{\mu}\subset  
\overline{\cal F}_{r+1} (f_{\mu})'_x  +  
\overline{\cal E}_{>r} (f_{\mu})'_y\ .$$
Plus pr\'ecis\'ement, si $\mbox{\em{in}}_r\,w= x^r W \in x^rE$, alors 
$\mbox{\em{in}}_{r+1}\,u = x^{r+1} (A_{\mu|F})^{-1}\,W$. 
\end{prop}

\begin{demo}
Nous conservons les notations du paragraphe \ref{parcalcwf} ; 
en particulier~:  $v= \sum_{i=1}^{n} v_i \varepsilon_i$, et
$w  = \sum_{i=1}^{n} w_i \varepsilon_i$ avec $\sum_{i=1}^{n} w_i = 0$.
Apr\`es division par $f_{\mu}\times \prod_{i=1}^n (y - a_i)$, 
l'\'equation (\ref{eqwfmult}) devient~:
$$\sum_{i,j}\; \frac{-u_i \mu_j a'_j + v_i \mu_j}{(y-a_i)(y-a_j)} \; \; = \; \; \sum_{i=1}^n \; \frac{w_i}{y-a_i}\ .$$
\`A nouveau, nous avons~: $v_i = u_i a'_i$, $1\leq i\leq n$, et :
$$\sum_{i\neq j} \frac{-u_i \mu_j (a'_i - a'_j)}{a_i - a_j} \left[ 
\frac{1}{y-a_i}  - \frac{1}{y-a_j}\right] = 
\sum_{i=1}^n \frac{w_i}{y-a_i}\ .$$
Par suite : 
$$ w_i = ( \sum_{j\neq i}\; \mu_j \frac{a'_i - a'_j}{a_i - a_j})  u_i  
- \sum_{j\neq i}( \mu_i \frac{a'_i - a'_j}{a_i - a_j})  u_j\ ,\ \, 
1\leq i\leq n\ .$$
Si $\cal U$ et $\cal W$ sont les matrices colonnes des
coordonn\'ees de $u$ et $w$, cette identit\'e s'\'ecrit :
${\cal W}  =  {\cal A}_\mu\cdot  {\cal U}$. De la 
m\^eme mani\`ere que lors de la preuve de la proposition
\ref{ArestFauto}, nous avons la d\'ecomposition~:
$ {\cal A}_\mu= (1/x)A_{\mu} +  {\cal A'}$, 
où $A_{\mu}$ est la matrice magique g\'en\'eralis\'ee 
associ\'ee au bouquet de branches avec multiplicit\'es
 d\'efini par $f_{\mu}$ (D\'efinition \ref{matmagen}),
et ${\cal A'}$ est une matrice \`a coefficients dans $\overline{K}$ 
de valuation strictement sup\'erieure \`a $-1$. Comme la
restriction de $A_\mu$ \`a $F$ est inversible (Proposition 
\ref{ValpropAmu}), la matrice ${\cal A}_{\mu}$ 
induit un endomorphisme inversible de $\overline{\cal F}$ 
dont le gradu\'e $\mbox{gr}\,{\cal (A_{\mu|\overline F})}^{-1}$ 
s'identifie \`a   $x  (A_{\mu|F})^{-1}$.
\end{demo}

\begin{coro} 
 Soit $f\in {\bf C}[[x]][y]$ un polynôme distingu\'e 
de degr\'e $n\geq 2$.
Notons $d\in {\bf N}$ le p.p.c.m des degr\'es de ses facteurs irr\'eductibles. 
 Pour tout rationnel  $r \in (1/d){\bf Z}$ :
$${\cal F}_{r}f \subset  {\cal F}_{r+1} f'_x  +  {\cal E}_{>r} f'_y\ .$$
\end{coro}

\subsection{L'op\'erateur $\nabla$}

Ce paragraphe est consacr\'e \`a l'\'etude de l'op\'erateur 
$\nabla:\overline{\cal F}\rightarrow\overline{\cal F}$ que nous
allons d\'efinir maintenant. 
Son utilit\'e sera patente dans la partie \ref{partbern}.

\begin{defi}
Soit  $f= \prod_{i=1}^{n} (y - a_{i})\in\overline{K}[y]$ un
polyn\^ome r\'eduit d\'efinissant un bouquet de branches. 
\`A tout $w\in\overline{\cal F}$, on associe 
$\nabla\,w\in\overline{\cal F}$ d\'efini par~: 
$$\nabla w\;=\; - \; u'_x \, - v'_y $$
o\`u $(u,v)$ est la solution dans 
$ \overline{\cal F} \times \overline{\cal E}$ de 
l'\'equation $ u f'_x + v f'_y  = w f $.
\end{defi}

Les r\'esultats obtenus au paragraphe pr\'ec\'edent
vont bien \'evidemment nous permettre de pr\'eciser
$\mbox{gr}\,\nabla$. Faisons d'abord une remarque
pr\'eliminaire.

\begin{lemm} \label{petitlemm}
Soit  $f= \prod_{i=1}^{n} (y - a_{i})\in\overline{K}[y]$ un
polyn\^ome r\'eduit d\'efinissant un bouquet de branches.
Soit $w\in\overline{\cal F}$ et $(u,v)$ la solution dans 
$\overline{\cal F} \times \overline{\cal E}$ de 
l'\'equation $ u f'_x + v f'_y  = w f $. Alors~:
$$\nabla\, w \; = \; - \, \left[\sum_{i=1}^n \, u_i' \, \varepsilon_i \right]  - \,w$$
o\`u $\sum_{i=1}^nu_i\varepsilon_i$ est la d\'ecomposition de $u$
dans la base $(\varepsilon_1,\ldots,\varepsilon_n)$.
\end{lemm}

\begin{demo}
En utilisant que 
$f = (y-a_i) \varepsilon_i$, $1\leq i\leq n$, nous obtenons 
les identit\'es :
\begin{equation} \label{derivlog}
\frac{f'_x}{f} = -\frac{a'_i}{y-a_i}+ \frac{(\varepsilon_i)'_x}
{\varepsilon_i} \; \ ; \ \; 
\frac{f'_y}{f} =  \frac{1}{y-a_i}+\frac{(\varepsilon_i)'_y}{\varepsilon_i}\ .
\end{equation}
Par ailleurs,  nous avons~: 
$v=\sum_{i=1}^n a'_i u_i\varepsilon_i$ 
(Lemme \ref{lemmprepawf}). Il
vient alors~:
 $u_i (\varepsilon_i)'_x  +  v_i(\varepsilon_i)'_y   = 
 (u_i f'_x  +  v_i f'_y)(\varepsilon_i/f)$, et donc :
$\sum_{i=1}^n  (u_i (\varepsilon_i)'_x  +  v_i(\varepsilon_i)'_y) \; = \; 
(u f'_x \, + \, v f'_y )/f \; = \; w .$ L'assertion en r\'esulte sans peine.
\end{demo}

\begin{prop} \label{grnabla}
  Soit $f=\prod_{i=1}^n(y-a_i)$ un polyn\^ome d\'efinissant un  
bouquet passant par l'origine et 
$A$ la matrice magique associ\'ee. Pour tout rationnel $r\in{\bf Q}$,
  l'op\'erateur $\nabla$  induit un endomorphisme 
de ${\overline {\cal F}}_r$. De plus~:
$$\mbox{\em{gr}}_r \nabla  \; = \; - \,   I \, 
- \, (1+r) (A_{|F})^{-1} $$
sous l'identification de $\mbox{\em{gr}}_r \overline{\cal F}$ 
avec $x^rF$.
\end{prop}

\begin{demo}
Notons $\cal U$, $\cal W$, ${\cal U}'$,  et ${\cal W}'$, 
les matrices colonnes form\'ees respectivement par 
les coordonn\'ees de $u$, de $w$ et de leur d\'eriv\'ees partielles en $x$
dans la base $(\varepsilon_1,\ldots,\varepsilon_n)$ de 
$\overline{\cal E}$. D'apr\`es la preuve de la proposition
\ref{ArestFauto} et l'identit\'e~: ${\cal W}={\cal A}\cdot{\cal U}$ 
(Lemme \ref{lemmprepawf}), nous avons~: 
${\cal U} =  x\left((A_{|F})^{-1} +  {\cal S}\right)\cdot{\cal W}$
o\`u $A$ est la matrice magique associ\'ee au bouquet et ${\cal S}$ 
est une matrice \`a coefficients de valuation strictement positive.
En d\'erivant, il vient :
$${\cal U}' =\left((A_{|F})^{-1}+{\cal T}\right) \cdot {\cal W} \; + \;
 x \left((A_{|F})^{-1} \;+\;  {\cal S}\right) \cdot {\cal W}'$$
avec ${\cal T} = {\cal S}+ x{\cal S}'_x$ \`a coefficients de 
valuation strictement positive. En particulier, si $w$ appartient 
\`a ${\overline {\cal F}}_r$,  alors $\nabla w $ y appartient aussi.
De plus,  si  $U'$ et  $W$ d\'esignent les formes initiales en degr\'e $r$
des matrices colonnes  ${\cal U}'$ et $\cal W$, nous avons 
alors :
$ U' = (1+r) \, (A_{|F})^{-1} \; W$, et avec le lemme pr\'ec\'edent~:
$\mbox{gr}\,\nabla \,( W) \, = 
\, - \, \left( I + (1+r) \, (A_{|F})^{-1} \, \right)\cdot W $. 
\end{demo}

Int\'eressons-nous maintenant au cas 
d'un bouquet de branches avec multiplicit\'es
 passant par l'origine, d\'efini par  
$ f_\mu =  \prod_{i=1}^n (y - a_i)^{\mu_i}$.
D'apr\`es ce qui a \'et\'e expliqu\'e au paragraphe \ref{resolwfnonred}, 
$\nabla$ est encore bien d\'efini comme endomorphisme
de ${\overline {\cal F}}$.

\stl

Les r\'esultats de la proposition pr\'ec\'edente s'\'etendent
sans peine.

\begin{prop}
  Soit $f_\mu=\prod_{i=1}^n(y-a_i)^{\mu_i}$ 
un polyn\^ome d\'efinissant un bouquet de branches avec
 multiplicit\'es, passant par l'origine. Notons  $A_{\mu}$ 
la matrice magique g\'en\'eralis\'ee associ\'ee et $A$ 
la matrice magique associ\'ee au bouquet de $n$ branches
sous jacent. Pour tout rationnel $r\in{\bf Q}$, l'op\'erateur 
$\nabla$ induit un endomorphisme de 
${\overline {\cal F}}_{r}$. De plus~: 
$$\mbox{\em{gr}}_r \nabla  \; = \; - \, \left( \, A \, + \, (1+r)I \, \right) \, (A_{\mu|F})^{-1}$$
sous l'identification de $\mbox{\em{gr}}_r \overline{\cal F}$ 
avec $x^rF$.
\end{prop}

\begin{demo} %
  Nous reprenons les notations utilis\'ees lors
de la preuve de la proposition \ref{grnabla} ; posons encore
$\overline{w}=
\sum_{i=1}^{n} (u_i (\varepsilon_i)'_x + v_i(\varepsilon_i)'_y)\in
\overline{\cal F}$ et $f=\prod_{i=1}^n(y-a_i)$.
Le calcul men\'e au lemme \ref{petitlemm} conduit \`a :
$ \overline{w}= (uf'_x +  v f'_y )/f$, puis :
\begin{equation} \label{nablawbar}
 \nabla w = -\left[\sum_{i=1}^n u'_i\varepsilon_i\right]-\overline{w}\ .
\end{equation}
En particulier, la solution 
$(u,v)\in\overline{\cal F} \times \overline{\cal E}$ 
de l'\'equation $ u(f_\mu)'_x  +  v (f_\mu)'_y  = w f_\mu$
 est aussi la solution de 
$ u f'_x + v f'_y  =  \overline{w} f$. 
D'apr\`es le lemme \ref{lemmprepawf}, nous avons donc :
$\overline{\cal W }=  (1/x) \left( A+{\cal A}'\right) \cdot {\cal U}$, 
o\`u $\overline{\cal W}$ est la matrice colonne des coordonn\'ees 
de $\overline{w}$, et ${\cal A}'$ est une matrice 
\`a coefficients dans $\overline{K}$ de valuation strictement positive. 
Par ailleurs, d'apr\`es la preuve de la proposition 
\ref{resolwfmultiple}, nous avons aussi :
${\cal U} =  x \left((A_{\mu|F})^{-1}+{\cal S}\right) \cdot {\cal W}$, et
donc encore~:
$${\cal U}' \; =  \left((A_{\mu|F})^{-1} \;+\;{\cal T}\right)\cdot{\cal W} 
+ x \, \left((A_{\mu|F})^{-1} \;+\;  {\cal S}\right) \cdot {\cal W}'$$
o\`u ${\cal S}$ et ${\cal T}$ sont des matrices \`a coefficients de 
valuation strictement positive.

\stl

  \`A partir de l'identit\'e (\ref{nablawbar}), nous en d\'eduisons 
alors que le vecteur colonne  des coordonn\'ees de $\nabla\,w$
est donn\'e par la somme~:
$$- \, \left( \, (A +I) \, (A_{\mu|F})^{-1} 
\, + \, {\cal L} \right)\cdot {\cal W}\; - \; x
\left((A_{\mu|F})^{-1} \;+\;  {\cal S}\right) \cdot {\cal W}'$$ 
où ${\cal L}$ est une matrice \`a coefficients de valuation strictement positive. En cons\'equence, $\nabla w$ appartient bien \`a 
$\overline{\cal F}_r$ lorsque $w\in{\overline {\cal F}}_r$,
et sa forme initiale  en degr\'e $r$ est :
$ - \, \left( \, A \, + \, (1+r)I \, \right) \, (A_{\mu|F})^{-1} \cdot W $,
o\`u $W=\mbox{in}_r\,w\in x^rF$. 
\end{demo}

\subsection{La r\'esolution de l'\'equation $u f'_x+v f'_y = w$}

Soit $f= \prod_{i=1}^n (y - a_{i})\in\overline{K}[y]$ un polyn\^ome 
r\'eduit. Dans ce paragraphe, nous nous proposons d'\'etudier 
l'\'equation :
\begin{equation}
  u f'_x+ \, v f'_y = w 
\end{equation}
pour $w$ donn\'e dans $\overline{\cal E}$. 
Nous savons que cette \'equation admet une unique 
solution $(u,v)$ dans $\overline{\cal F}\times \overline{\cal F}$. 
Avant de la d\'eterminer pr\'ecis\'ement  dans le cas d'un bouquet 
passant par l'origine, fixons quelques notations.

\begin{nota} {\em
  \'Etant donn\'e un bouquet \`a $n$ branches 
$f=\prod_{i=1}^n (y-a_i)$ 
passant par l'origine,  nous notons : 
$\pi = \mbox{ sup} \{ m_i + m_{i,j} \, ; \, 1\leq  i,\,j\leq n,\ i\not=j\}$, 
$\delta = \mbox{ inf} \{ \nu(a_i) \, ; \, 1\leq  i\leq n\}$,
$\theta = \pi \, - \, 1$, et $\tau \, = \, \pi \, - \delta$.}
\end{nota}

\begin{prop}  \label{laresolw} 
  Soit $f=\prod^n_{i=1}(y-a_i)\in\overline{K}[y] $ un 
polyn\^ome r\'eduit d\'efinissant un bouquet 
passant par l'origine.  

(i) Pour tout $w$ appartenant \`a $\overline{\cal E}$, 
la solution dans $\overline{\cal F} \times \overline{\cal F}$ 
de l'\'equation $ u f'_x +  v f'_y = w $ est donn\'ee par :
$$\left\{  \begin{array}{ccl}
    u  & = &  {\cal (A_{| \overline F})}^{-1} {\cal B} w   \\ 
    v  & = & \sum_{i=1}^n (u_i a'_i  +(w_i/\varepsilon_i(a_i)) ) 
  \varepsilon_i \end{array}\right. $$
o\`u $u=\sum_{i=1}^n u_i \varepsilon_i$, 
$w=\sum_{i=1}^n w_i\varepsilon_i$, et 
$\cal B$ d\'esigne l'endomorphisme de $\overline{\cal E}$ induit par
la matrice magique de ${\cal M}ag_0(n,\overline{K})$, sym\'etrique,  
de terme  g\'en\'eral~:
$\beta_{i,j}=-(\varepsilon_i(a_i) - \varepsilon_j(a_j))/
((a_i - a_j)\cdot \varepsilon_i(a_i) \cdot \varepsilon_j(a_j))$, $i\not=j$.

En particulier, pour tout $r\in {\bf Q}$ :
$$\overline{\cal E}_{r} \; \subset \; \overline{\cal F}_{r-\theta} \, f'_x \; + \; \overline{\cal F}_{r-\tau} \, f'_y\ .$$
(ii) Soit $w=uf'_x+vf'_y$ un \'el\'ement de $\overline{\cal E}_r$
avec $u\in \overline{\cal F}_{r-\theta}$ et 
$v\in \overline{\cal F}_{r-\tau}$. Alors $u'_x+v'_y$ appartient
\`a $\overline{\cal F}_{r-\pi}$.
\end{prop}

\begin{demo}
 Posons  $v =\sum_{i=1}^n v_i \varepsilon_i$ avec 
$\sum_{i=1}^n v_i = 0$. En divisant l'\'equation par $f^2$ puis
 en utilisant l'identit\'e~: 
$1/f=\sum_{j=1}^n (1/\varepsilon_j(a_j))(1/(y-a_j))$, 
 nous obtenons  :
$$\sum_{i,j}\; \frac{-u_i a'_j + v_i}{(y-a_i)(y-a_j)} \; \; = \; \; \frac{1}{f} \, \left(\sum_{i=1}^n \; \frac{w_i}{y-a_i}\right)
\; \; = \sum_{i,j}\; \frac{w_i}{\varepsilon_j(a_j)} \, \cdot \, \frac{1}{(y-a_i)(y-a_j)}$$
Par unicit\'e de la d\'ecomposition en \'el\'ements simples, 
l'identification fournit :
\begin{equation} \label{vfoncu}
\frac{w_i}{\varepsilon_i(a_i)} =  v_i - u_i a'_i\ ,\ \, 1\leq i\leq n\ \; .
\end{equation}
Apr\`es \'elimination  des $v_i$, l'\'equation devient :
$$\sum_{i\neq j}\; \frac{u_i(a'_i - a'_j)}{(y-a_i)(y-a_j)} \; \; = \; \; \sum_{i \neq j} \; w_i
(\frac{1}{\varepsilon_j(a_j)} - \frac{1}{\varepsilon_i(a_i)}) \cdot \frac{1}{(y-a_i)(y-a_j)}\ .$$
En utilisant l'identit\'e (\ref{fracrat}) page 
\pageref{fracrat}, il vient finalement :  
$$\sum_{i\neq j}\; (\frac{a'_i - a'_j}{a_i - a_j}) \, u_i \; \cdot \,\left[\frac{1}{y-a_i} -\frac{1}{y-a_j}\right] \;\; = 
\mbox{\ \ \ \ \ \ \ \ \ \ \ \ \ \ \ \ \  }$$
$$\mbox{\ \ \ \ \ \ \ \ \ \ } \sum_{i \neq j} \; \frac{\varepsilon_i(a_i) - \varepsilon_j(a_j)}{a_i - a_j}\, \cdot \, \frac{w_i}{\varepsilon_i(a_i) \varepsilon_j(a_j)} \, \cdot \, \left[\frac{1}{y-a_i} -\frac{1}{y-a_j}\right]\ .$$
Si $\cal U$ et $\cal W$ d\'esignent encore 
les matrices colonnes form\'ees par les coordonn\'ees de $u$ et $w$
dans la base $(\varepsilon_1,\ldots, \varepsilon_n)$ de 
$\overline{\cal E}$, cela se r\'e\'ecrit~:
$ \cal A \, \cal U \; = \;\cal B \,\cal W $,
où $\cal A$ est la matrice magique compl\`ete 
associ\'ee au bouquet. Par ailleurs, nous remarquons
que l'image de $\cal B$ est contenue dans $\cal F$.
Gr\^ace \`a la proposition \ref{ArestFauto},
nous obtenons alors l'expression de $(u,v)$ annonc\'ee.

\stl

Observons maintenant que $\beta_{i,j}$ est de valuation sup\'erieure 
ou \'egale \`a $- \pi$. Comme la valuation des termes  
de la matrice de l'endomorphisme ${\cal (A_{| \overline F})}^{-1}$  
dans la base $(\varepsilon_1-\varepsilon_n,\ldots,
\varepsilon_{n-1}-\varepsilon_n)$  est au moins \'egale \`a $1$ 
(Proposition \ref{ArestFauto}), 
nous en d\'eduisons que : $\mbox{val}(u) \geq \mbox{val}(w) - \theta$ 
et $\mbox{val}(v) \geq \mbox{val}(w)  -\tau $ ; d'o\`u la seconde  
assertion.

\stl

Montrons maintenant le second point. 
Compte-tenu  des identit\'es (\ref{derivlog}) et  (\ref{vfoncu}),
nous avons~:
$$\sum_{i=1}^{n} \; (u_i (\varepsilon_i)'_x \, + \, v_i(\varepsilon_i)'_y) \; = \; \frac{w}{f} \; - \; \sum_{i=1}^{n} 
\frac{w_i}{\varepsilon_i(a_i)} \frac{\varepsilon_i}{y-a_i}\ .$$
Par suite :
$$ u'_x+v'_y\; = \; \sum_{i=1}^n \, \, u_i' \, \varepsilon_i \;  
- \, \sum_{i=1}^{n} 
\frac{w_i}{\varepsilon_i(a_i)} \frac{\varepsilon_i - 
\varepsilon_i(a_i)}{y-a_i}.$$
Remarquons alors que, 
 pour des indices $i$ et $j$ distincts, 
  la valuation de $((\varepsilon_i - \varepsilon_i(a_i))/(y-a_i))(a_j) = -\varepsilon_i(a_i)/(a_j-a_i)$ est \'egale \`a $m_i-m_{i,j}$ ; quant \`a
celle de 
$((\varepsilon_i - \varepsilon_i(a_i))(y-a_i))(a_i) = (\varepsilon_i)'_y(a_i)$,
elle est  sup\'erieure ou \'egale \`a $m_i - \mbox{sup}_{k \neq i}\{m_{i,k}\}$. 
Avec la remarque \ref{Ezero},  nous en d\'eduisons alors  que :
$\mbox{val}((\varepsilon_i - \varepsilon_i(a_i))/(y-a_i))  \geq  
\mbox{inf}\{ \, \mbox{inf}_{j\neq i} \{ m_i - m_{i,j} - m_j \} \, , \, \mbox{inf}_{k\neq i} \{\,- m_{i,k}\}\}$. 
D'où : $\mbox{val}((w_i/\varepsilon_i(a_i))(\varepsilon_i - \varepsilon_i(a_i))/(y-a_i))  \geq \nu(w_i) 
- \mbox{sup}_{j \neq i} \{ m_{i,j} + m_j \}$, puis~: 
$\mbox{val}(u'_x+v'_y) \geq \mbox{inf} \{ \mbox{val}(u)-1 \, , \, \mbox{val}(w)-\pi \}$. La partie $(ii)$ s'ensuit.
\end{demo}

Le corollaire \ref{cororesolwf} a  son pendant.

\begin{coro} \label{cororesolw} \label{coronabla}
  Soit $f\in{\bf C}[[x]][y]$ un polynôme distingu\'e r\'eduit. 

(i) Pour tout $r\in {\bf Q}$ : 
$${\cal E}_{r} \; \subset \; {\cal F}_{r-\theta} \, f'_x \; + \; {\cal F}_{r-\tau} \, f'_y\ .$$
En particulier, ${\cal E}_{\pi}$ est contenu dans l'id\'eal jacobien 
$(f'_x,f'_y){\bf C}[[x]][y]$.

(ii) Soit $w=uf'_x+vf'_y$ un \'el\'ement de ${\cal E}_\pi$
avec $u\in {\cal F}_{1}$ et 
$v\in {\cal F}_{\delta}$. Alors $u'_x+v'_y$ appartient
\`a ${\cal F}_{0}$.
\end{coro}

\begin{rema}{\em
Si l'on souhaite \'etendre ces r\'esultats au cas d'un 
polynôme $g$ de degr\'e sup\'erieur ou \'egal \`a $n$,
il convient d'abord de le d\'ecomposer selon les puissances 
successives de $f$~: $\, 
g \, = \, w_0 + w_1 f + \cdots + w_{\ell} f^{\ell}$ avec 
$w_i\in\overline{\cal E}$, $0\leq i\leq \ell$. Il s'agit ensuite de
 trouver une minoration convenable des poids des $w_k $, 
pour $1\leq k \leq \ell$ :
$\mbox{val}(w_k) \, = \, \mbox{inf} \{\nu(w_k(a_i)/\varepsilon_i(a_i))
 \; ; \; 1 \leq i \leq n \}$ en fonction de la donn\'ee $g$. 
Par exemple, si $g$ est un polyn\^ome de degr\'e strictement 
inf\'erieur \`a $2n$,  alors : $g=w_0 + w_1 f$, et
nous pouvons minorer les poids de $w_0$ et $w_1$ 
\`a l'aide des valuations des $g(a_i)$ et $g'_y(a_i)$, $1 \leq i \leq n$, 
et des multiplicit\'es d'intersection des branches du bouquet, en utilisant 
les identit\'es :
$$g(a_i) \; = \; w_0(a_i) \; \;   ; \ \; g'_y(a_i) \; = \; 
(w_0)'_y(a_i) \, + \, w_1(a_i) \, f'_y(a_i)\ .$$
Apr\`es calcul, nous trouvons :
$\mbox{val}(w_0) \, = \mbox{inf} \, \{ \nu(g(a_i))-m_i \; ; \; 1 \leq i\leq n \}$, 
et le poids   $\mbox{val}(w_1)$ est sup\'erieur ou \'egal \`a~:
$$\mbox{inf}\{\nu(g'_y(a_i))-2m_i, \, \nu(g(a_j))-m_i-m_j-m_{i,k} \; ; \; 1 \leq i,j,k\leq n,\; i\not=k \}\ .$$ 
On peut alors appliquer la proposition \ref{laresolw} \`a $w_0 $ et 
le corollaire \ref{laresolwf} \`a $w_1 f$.}
\end{rema}

\begin{exem}{\em 
\`A partir du corollaire \ref{cororesolw} et de la remarque pr\'ec\'edente, 
nous sommes en mesure de pr\'eciser  une puissance de
 l'id\'eal maximal contenue dans l'id\'eal jacobien d'un germe de 
courbe \`a singularit\'e isol\'ee $f $
en fonction des multiplicit\'es des composantes irr\'eductibles et des multiplicit\'es d'intersection des composantes entre
elles. Cette puissance ne d\'epend donc que du type topologique du germe, et elle est meilleure\footnote{C'est-\`a-dire plus petite} que le nombre de Milnor,  
$\mu=\sum_{i=1}^n m_i-n+1$, qui convient toujours.}
\end{exem}

\newpage

\section{Un multiple de la $b$-fonction}
\label{partbern}

\subsection{Pr\'eliminaires}

\subsubsection{Le polynôme de Bernstein}

 Notons ${\cal O}$ l'anneau des germes de fonctions
 holomorphes \`a l'origine de ${\bf C}^2 , \; {\cal D}$ l'anneau des
germes d'op\'erateurs diff\'erentiels \`a coefficients dans ${\cal O}$ et ${\cal D}[s] \; = \; {\cal D} \otimes _{\bf C} {\bf C} [s] .$

\stl

\'Etant donn\'e un germe $f$ de fonction holomorphe, nous consid\'erons le ${\cal O} [1/f , s]$-module libre de rang un, 
${\cal O} [1/f, s] f^s$, 
muni de sa structure naturelle de ${\cal D}[s]$-module. 
D'apr\`es \cite{K}, il existe un polynôme non nul $b(s)\in{\bf C}[s]$ 
et un op\'erateur $P\in{\cal D}[s]$ r\'ealisant dans 
${\cal O} [1/f, s] f^s$ l'\'equation fonctionnelle :
$$ b(s) f^s \; = \; P\cdot f^{s+1}\ .$$
Le polyn\^ome unitaire $b(s)$ de plus petit degr\'e 
r\'ealisant une telle 
\'equation est appel\'e {\em polynôme de Bernstein} de $f$ 
ou {\em $b$-fonction} de $f$. Nous renvoyons le lecteur 
aux travaux de B. Malgrange (\cite{Mal1}, \cite{Mal2}), 
A. N. Varchenko (\cite{Var}) et T. Yano (\cite{Ya}) sur ce sujet.

\stl

Nous supposons maintenant que $f$ s'annule \`a l'origine et 
d\'efinit un germe de singularit\'e isol\'ee. Son polynôme
de Bernstein $b(s)$ est alors de la forme 
$b(s)=(s+1)\tilde{b}(s)$ où $\tilde{b}(s)$ est le
polynôme minimal de l'action de $s$ sur le 
${\cal D}$-module de type fini :
$$ {\cal M} \; = \; (s+1) \frac{{\cal D}[s] f^s}{{\cal D}[s] f^{s+1}} \; \cong
 \; \frac{{\cal D}[s] f^s}{{\cal D}[s](f,f'_x,f'_y)f^s}$$
ou sur sa cohomologie de de Rham :
$$ H^2 {\cal M} \; = \; \frac{{\cal M}}{\partial_x {\cal M} + \partial_y {\cal M}}$$
qui est alors un ${\bf C}$-espace vectoriel de dimension $\mu$, le nombre de Milnor de la singularit\'e (voir \cite{Mal1}).

\stl

  Nous  supposerons de plus que $f$ est un polynôme distingu\'e 
r\'eduit de ${\bf C}\{x\}[y]$ de degr\'e $n\geq 2$. D'apr\`es le 
th\'eor\`eme de Newton-Puiseux, il s'\'ecrit :
$f= \prod_{i=1}^{n} (y - a_{i})$ dans ${\bf C}\{x\}[x^{1/d}][y]$
 pour l'entier $d$ \'egal au p.p.c.m.
des degr\'es des facteurs irr\'eductibles de $f$. 
Dans cette partie, nous allons construire un multiple
du polynôme de Bernstein de $f$ qui ne d\'epend que de 
$d$ et de la famille $\{ m_{i,j} \; ; \; 1\leq i < j\leq n \}$
des multiplicit\'es d'intersection des branches deux-\`a-deux. 
Lorsque $n$ est la multiplicit\'e\footnote{C'est-\`a-dire lorsque l'axe 
des $y$ est transverse \`a $f$.} de $f$, ce 
multiple ne d\'epend alors que du type topologique
du germe de courbe d\'efini par $f$. 

\stl

Dans cette construction, nous n'utiliserons pas 
la cohomologie de de Rham invoqu\'ee ci-dessus. Par contre, nous 
verrons que cette cohomologie $ H^2 {\cal M}$ est engendr\'ee 
par les classes d'\'el\'ements de la forme
$w f^s$ et $s w f^{s-1}$ avec $w$ dans ${\cal O}$.

\subsubsection{La transformation de Tchernhaus}

Soit $f\in{\bf C}\{x\}[y]$ un polyn\^ome distingu\'e, 
r\'eduit, et $\prod_{i=1}^n(y-a_i)$
sa d\'ecomposition dans ${\bf C}\{x\}[x^{1/d},y]$.
Nous rappelons que $\sigma$ d\'esigne la somme 
$\sum_{i=1}^n a_i$. Posons :
 $$ \alpha = \mbox{inf} \{ \, m_{i,j} \; \; ; \; \; 1\leq i < j\leq n \}\ .$$ 
Nous constatons que le rationnel  $\alpha$ 
co\"{\i}ncide avec 
$\mbox{inf}\{\nu(a_i -\sigma/n)\; ; \; 1\leq i \leq n \}$. C'est 
une cons\'equence des deux in\'egalit\'es suivantes :

- pour $1\leq i\leq n$~: $\nu(a_i - \sigma/n) \geq \alpha$ 
 puisque $a_i - \sigma/n \; = \; \sum_{j=1}^n(a_i - a_j)/n$ ;

- pour $1\leq i\not=j\leq n$~: $\nu(a_i - a_j)  \geq \mbox{inf}\{ \nu(a_i - \sigma/n)\; ; \; 1\leq i \leq n \}$ ayant $\; a_i - a_j = (a_i - \sigma/n) - (a_j - \sigma/n)$.

\stl

Nous en d\'eduisons que le polynôme $h\in{\cal F}$ d\'efini par~:
\begin{eqnarray*}
   h & = & f \, - \, (\frac{y}{n}-\frac{\sigma}{n^2}+\frac{x\sigma'}{n^2\alpha}) \, f'_y \, - \,  \frac{x}{n\alpha} \, f'_x \\
      & = & \sum_{i=1}^n \frac{-1}{n\alpha}\left[\alpha (a_i - \frac{\sigma}{n})  - x (a'_i - \frac{\sigma'}{n})\right]\varepsilon_i 
\end{eqnarray*}
est de poids $\mbox{val}(h)$ strictement sup\'erieur \`a $\alpha$. 
Nous faisons alors le changement de coordonn\'ees : 
$x_1 = x \, , \, y_1 = y - {\sigma}/{n}$. Pour ne pas alourdir les
notations, nous appelons toujours $x$ et $y$ ces nouvelles coordonn\'ees et nous supposerons d\'esormais que $\sigma$ est nul. 
En particulier, nous avons alors :
 $$ \alpha = \mbox{inf}\{ \, m_{i,j} \; ; \; 1\leq i < j\leq n \} = 
  \mbox{inf} \{ \, \nu(a_i) \; ; \; 1\leq i \leq n \} = \delta\ , $$
$$ h \; = \; f \, - \, \frac{y}{n} \, f'_y \, - \, \frac{x}{n\delta} \, f'_x \;
= \;\frac{1}{n\delta}\sum_{i=1}^n(x a'_i -  \delta a_i ) \varepsilon_i$$
avec $\mbox{val}(h) > \delta$. 

\stl

Pour $1\leq i\leq n$, notons $a_{i,\delta}\in{\bf C}$ le coefficient de 
$x^{\delta}$ dans le d\'eveloppement de $a_i$ ; alors :
 $\nu(a_i - a_{i,\delta} x^{\delta}) > \delta $. Constatons enfin que
le premier côt\'e du polygone de Newton de $f$ a pour \'equation :
 $k/n\delta + l/n = 1$ dans les coordonn\'ees $(k,l)$, et que 
 la restriction de $f$ \`a ce côt\'e est le polynôme quasi-homog\`ene : 
$\prod_{i=1}^{n} (y - a_{i,\delta} x^{\delta})$.

\subsubsection{La division selon le premier c\^ot\'e du polygone de
Newton} 
 
Introduisons maintenant la fonction de poids 
$\rho:{\cal O}\rightarrow {\bf Q}^+$ associ\'ee au
premier c\^ot\'e du polygone de Newton de $f$.
Pour tout $w=\sum w_{k,l}x^ky^l\in{\cal O}$
non nul, nous posons~:

$$ \rho (w) \; = \; \mbox{inf}\, \{ \; \frac{k}{n\delta} + \frac{l}{n} \; ; \; w_{k,l} \neq 0 \}\ .$$
La filtration de ${\cal O}$ associ\'ee est d\'efinie
alors en posant,  pour tout $q\in {\bf Q}^+$~:
${\cal O}_{q}=\{w\in{\cal O}\, ; \, \rho(w)\geq q \}$ - avec la
convention~: $\rho(0)=+\infty$.
Notons aussi ${\cal O}_{>q}=\{w\in{\cal O}\, ; \, \rho(w)>q \}$. 

\stl

Nous avons bien s\^ur~:
$\rho (f) = 1$. D'autre part, la relation suivante permet de
comparer les poids  $\rho$ et $\mbox{val}$ d'un \'el\'ement
de ${\cal E}\cap{\cal O}=\sum_{i=0}^{n-1}{\bf C}\{x\}y^i$.

\begin{lemm} \label{comparpoids}
  Pour tout $w\in {\cal E}\cap {\cal O}$~:
$$\mbox{\em val}(w)\geq n\delta \rho(w)-\mbox{\em sup}\{m_i\, ;\,
1\leq i\leq n \}\  . $$
En particulier : ${\cal E}\cap{\cal O}_q  \subset {\cal E}_0$ d\`es que
 $q\geq{\mbox{\em sup}\{m_i \}/n\delta}$.
\end{lemm}

\begin{demo}
  Comme $\nu(x^k a_{i}^l) = n\delta \rho (x^k y^l)$ 
pour tout indice $i$ tel que $\nu(a_i) = \delta$, nous avons~: 
$\nu(w(x,a_i)) \geq n\delta \rho (w)$ pour $1\leq i\leq n$.
Ayant de plus $\mbox{val}(w)  =  \mbox{inf} \{\nu(w(x,a_i)) - m_i \, ;  1 \leq i \leq n  \}$ (voir la remarque \ref{Ezero}), le r\'esultat s'ensuit.
\end{demo}

Pour tout $ w \in {\cal O}$, nous avons l'identit\'e :
\begin{equation} \label{montpoidsw}
\frac{x}{n\delta} \, w'_x \, + \,  \frac{y}{n} \, w'_y \; = 
\; \rho(w) w \, + \, w_1 \  , \ \; \rho(w_1)>\rho(w) \ .
\end{equation}
En particulier, pour $w = f$~:
\begin{equation} \label{montpoidsf}
\frac{x}{n\delta} \, f'_x \, + \,  \frac{y}{n} \, f'_y \; = \; f \, - \, h \  , \ \; \rho(h)>1 \ .
\end{equation}
Enfin, le r\'esultat suivant s'obtient tr\`es facilement, en suivant 
pas-\`a-pas l'algorithme de la division.

\begin{lemm} \label{divpoireau} 
Soit $w \in {\bf C}\{x\}[y]$ un polynôme non nul. 
La division euclidienne de $wh$ par $f$ fournit :
\begin{equation} \label{eqdivpoireau}
 w h \; = \; w_2 f \; + \; \lambda(x) f'_y \; + \; \tilde{w}
\end{equation}
o\`u $w_2\in{\bf C}\{x\}[y]$, $\lambda(x)\in{\bf C}\{x\}$ et 
$\tilde{w}\in{\cal F}$ v\'erifient :
$\nu(\lambda(x)) > n \delta \rho(w) + \delta$,
$\rho(w_2) \geq \rho(w) + \rho(h) - 1 > \rho(w)$,  et
$\mbox{\em val}(\tilde{w}) > n \delta \rho(w) +\delta$.
\end{lemm}

\subsection{La mont\'ee des poids}

Dans ce paragraphe, nous donnons les lemmes techniques
sur lesquels s'appuit la d\'etermination d'un multiple
du polyn\^ome de Bernstein de $f$. La m\'ethode
utilis\'ee - la  `mont\'ee des poids'  - est classique
pour  traiter ce probl\`eme (voir \cite{BGMM} par exemple).

\stl

Nous consid\'erons trois classes d'\'el\'ements dans ${\cal O}[1/f,s]f^s$ :

\stl

\noindent\ \,classe A\,: les \'el\'ements de la forme  $w f^s$ avec 
$w \in {\cal O} \cap {\cal F}$ et $\mbox{val}(w) \geq 0$ ;

\noindent\ \,classe B\, : les \'el\'ements de la forme  $w f^s$ 
avec $w \in {\cal O} \cap {\cal F}$ et $\mbox{val}(w) < 0$ ;

\noindent\ \,classe C\, : les \'el\'ements de la forme  $s w f^{s-1}$ 
avec $w \in {\cal O} \cap {\cal F}$ et $\mbox{val}(w)>0$.

\stl

Comme ${\cal F}_0$ est contenu dans ${\cal O}$ (Remarque 
 \ref{Ezero}), la classe A (resp. C) est l'ensemble  ${\cal F}_{\geq 0}f^s$ 
(resp. $s{\cal F}_{>0}f^{s-1}$). En ce qui concerne un \'el\'ement 
$wf^s$ de la classe B, $wf$ n'appartient pas - en g\'en\'eral - \`a 
l'id\'eal jacobien ; pour ces \'el\'ements, nous ne pourrons faire 
que le premier pas de la division par cet id\'eal, selon le premier c\^ot\'e du 
polygone de Newton.

\subsubsection{Mont\'ee des poids en classe A}

\begin{lemm}\label{lemA}
  Soit $w \in {\cal F }_r$ avec $r\geq 0$. Alors :
$$M_r(s) \, w f^s \;\in \; {\cal D}[s] \, {\cal F }_{>r} f^s \; + \; 
{\cal D}[s] \, f'_y f^s $$
avec $M_r(s)=\prod_{\lambda\in\Lambda(r)}(s-\lambda)$ o\`u
$\Lambda(r)$ d\'esigne l'ensemble des valeurs propres de
l'endomorphisme diagonalisable $\mbox{\em gr}_r\nabla$.
\end{lemm}

\begin{demo}
D'apr\`es le corollaire \ref{cororesolwf}, il existe $u\in{\cal F}_{r+1}$
et $v\in{\cal E}_{>r}$ tels que : $wf = u f'_x + v f'_y$. 
Par suite :
$$ s w f^s = u \frac{\partial}{\partial x} \, f^s + v \frac{\partial}{\partial y} \, f^s = \nabla(w) f^s +
\frac{\partial}{\partial x} \, u f^s +  \frac{\partial}{\partial y} \, v f^s .$$

Nous pouvons \'ecrire : $v = v_1 + \eta(x) f'_y $ avec 
$v_1\in {\cal F }_{>r}$ et $\eta(x) \in {\bf C} \{x\}$ tel
que $\nu(\eta(x))>r$. Pour tout nombre complexe $\lambda$, 
nous avons alors :
$$(s-\lambda) w f^s = (\nabla - \lambda I)(w) f^s + R$$
 avec $R \in {\cal D}[s]{\cal F}_{>r}f^s + {\cal D}[s] f'_y f^s $.
D'apr\`es la proposition \ref{grnabla}, 
nous savons que $\nabla(w)$ appartient \`a ${\cal F}_r$.
On conclut en it\'erant la formule pr\'ec\'edente lorsque 
$\lambda$ parcourt $\Lambda(r)$. 
\end{demo}
 
Nous rappelons  que, gr\^ace au corollaire 
\ref{AFcourbliss}, l'ensemble $\Lambda(r)$ 
est compl\`ete-ment d\'etermin\'e
par l'arbre associ\'e au bouquet d\'efini par $f$.

\subsubsection{Mont\'ee des poids en classe B}

\begin{lemm} \label{lemB}
  Soit $w \in {\cal O}_q\cap {\cal F}$ avec $q\in{\bf Q}^+$. Alors :
$$(s+\frac{1}{n\delta}+\frac{1}{n}+q) \, w f^s \; \in \; 
{\cal D}[s] ({\cal O}_{>q}\cap {\cal F}) f^s \; + \; 
{\cal D}[s] \, f'_y f^s \; + 
\; {\cal D}[s] \, s {\cal F}_{>n \delta q + \delta } f^{s-1}\ .$$
\end{lemm}

\begin{demo}
Nous avons l'identit\'e~:
\begin{eqnarray*}
\left[\frac{1}{n\delta} \frac{\partial}{\partial x} x \, + \, \frac{1}{n} \frac{\partial}{\partial y} y \right] wf^s  & = &
 (s+\frac{1}{n\delta}+\frac{1}{n}+q)wf^s \;  \\
& & + \; sw\left(\frac{x}{n\delta} \, f'_x \, + \,  \frac{y}{n} \, f'_y \, - \, f\right) f^{s-1} \\
& & + \; \left(\frac{x}{n\delta} \, w'_x \, + \,  \frac{y}{n} \, w'_y \, -\, q w\right) f^s\ .
\end{eqnarray*}  
Soit encore  $yw=w_3+\eta(x)f'_y$ avec 
$w_3\in {\cal O}_{>q}\cap {\cal F}$,
la d\'ecomposition de $yw$ dans ${\cal F}\oplus {\bf C}\{x\}f'_y$. 
En utilisant les identit\'es (\ref{montpoidsw}), (\ref{montpoidsf}), 
et (\ref{eqdivpoireau}), nous trouvons alors~:
$$(s+\frac{1}{n\delta}+\frac{1}{n}+q)wf^s \; = \; s\tilde{w}f^{s-1}
\mbox{\ \ \ \ \ \ \ \ \ \ \ \ \ \ \  \ \ \ \ \ \ \ \ \ \ \ \ \ \ \ \ \ \ \ \ \ \ \ \ \ \ \ \ \ \ \ \ \ \ \ \ } \;$$
$$\mbox{\ \ \ \ \ \ \ \ \ \ \ \ \  \ \ \ \ \ \ \ \ } + \; \left[ \frac{1}{n\delta} \frac{\partial}{\partial x} xw  + \frac{1}{n} \frac{\partial}{\partial y} (w_3 + \eta(x) f'_y) -  w_1  +  sw_2 \, + 
\frac{\partial}{\partial y}\lambda(x)\right] f^s $$
avec~: $\rho(xw),\,\rho(w_1),\,\rho(w_2),\,\rho(w_3)>q$,  
$\rho(\lambda(x))= (1/n\delta)\nu(\lambda(x))>q$ et 
$\mbox{val}(\tilde{w}) > n \delta q + \delta$. D'o\`u le r\'esultat 
annonc\'e. 
\end{demo}

\subsubsection{Mont\'ee des poids en classe C} 

\begin{lemm} \label{lemC}
  Soit $w \in {\cal F }_r$ avec $r\geq 0$. Alors :
$$M_r(s-1) \; s w f^{s-1} \; \in \; {\cal D}[s] \; {\cal F }_{>r}s f^{s-1} \; + \; 
{\cal D}[s] ({\cal O}_{>r/n\delta}\cap {\cal F}) f^s $$
o\`u $M_r(s)$ d\'esigne le polyn\^ome minimal de 
$\mbox{\em gr}_r\nabla$.
\end{lemm}

\begin{demo} 
Nous proc\'edons de la même mani\`ere qu'au lemme \ref{lemA},
apr\`es d\'ecalage de $s$ en $s-1$. D'apr\`es le 
corollaire \ref{cororesolwf}, 
il existe $u \in {\cal F}_{r+1}$ et $v\in{\cal E}_{>r}$ tels que : 
$wf = u f'_x + v f'_y$ ; par suite :
$$ (s-1) w f^{s-1} = u \frac{\partial}{\partial x} \, f^{s-1} + v \frac{\partial}{\partial y} \, f^{s-1}
 = \nabla(w) f^{s-1} +
\frac{\partial}{\partial x} \,u f^{s-1}+\frac{\partial}{\partial y} \, v f^{s-1} .$$
Nous pouvons \'ecrire: $v = v_1 + \eta(x) f'_y$ avec 
$v_1\in{\cal F }_{>r} $ et $\nu(\eta(x))>r$. Pour
 tout nombre complexe $\lambda$, nous avons donc :
$$ (s-1-\lambda) s w f^{s-1} = s(\nabla - \lambda I)(w) f^{s-1} +
\left(\frac{\partial}{\partial y}\right)^2 \eta(x)f^s + R$$
avec $R \in {\cal D}[s] {\cal F }_{>r} sf^{s-1}$. Constatons que : 
$\rho(\eta(x)) = (1/n\delta)\nu(\eta(x)) > r/{n\delta}$. 
L'assertion s'obtient alors en it\'erant la formule pr\'ec\'edente,
 $\lambda$ parcourant l'ensemble des valeurs propres de 
$\mbox{gr}_{r}\nabla$.
\end{demo}

\subsection{R\'esultat des courses}

\subsubsection{Les polyn\^omes $M_{gros}^A(s)$ et $M_{fin}^A(s)$}

  Un multiple grossier convenant \`a tous les \'el\'ements 
$wf^s$ en classe A est le produit $M_{gros}^A(s)$ des polyn\^omes $M_r(s)$ 
lorsque $r$ parcourt $(1/d)\bf{N}$ avec $0 \leq r < \pi$.  
Par it\'eration du lemme \ref{lemA} puis application 
du corollaire \ref{cororesolw}, nous obtenons :
\begin{equation} \label{eqgrosA}
 M_{gros}^A(s) {\cal F}_0f^s \; \subset {\cal D}[s] (f'_x,f'_y) f^s\ .
\end{equation}
Ce multiple ne d\'epend que de $d$ et des multiplicit\'es d'intersection 
$m_{i,j}$.

\stl

 Un multiple plus fin, $M_{fin}^A(s)$, s'obtient en prenant 
le produit de ces mêmes polynômes $M_r(s)$ lorsque
$r$ parcourt l'ensemble : 
$$\Gamma_f^0 \; = \; \{ \,\mbox{val}(w) \; ; \; w \in {\cal F}_0 \mbox{ et } w \notin (f,f'_x,f'_y)\} \ .$$
En g\'en\'eral, ce multiple ne d\'epend pas que du type topologique 
du germe de courbe d\'efini par $f$ (par exemple, penser
\`a $f$ quasi-homog\`ene et semi-quasi-homog\`ene). 
Nous avons alors :
$$ M_{fin}^A(s) {\cal F}_{0}f^s \; \subset {\cal D}[s] (f,f'_x,f'_y) f^s\ . $$

\subsubsection{Les polyn\^omes $M_{gros}^B(s)$ et $M_{fin}^B(s)$}

Nous rappelons que ${\cal O}_q\cap {\cal F}\subset {\cal F}_0$ 
pour  $q\geq \mbox{sup}\{m_i \}/n\delta$ (Lemme \ref{comparpoids}).
Pour d\'efinir un multiple grossier, nous consid\`ererons alors tous 
les poireaux :  $ \rho = k/n\delta + l/n$ avec   
$0\leq \rho <\mbox{sup}\{m_i\}/n\delta$ et $(k,l) \in {\bf N}^2$. 
Nous notons donc  $M_{gros}^B(s)$ le produit des facteurs 
$ (s + (k+1)/n\delta + (l+1)/n)$ correspondants. 
On remarquera que la plus petite valeur propre non nulle de la
matrice magique associ\'ee \`a $f$ est 
$\lambda_0=n\delta=\gamma(T_0)n(T_0)$, et que  d'autre part,
$r=k+(l+1)\delta$ est la valuation du reste de la division
par $f$ du polyn\^ome $wyf'_y$ avec $w=x^k y^{l}$. On peut donc r\'ecrire
$M^B_{gros}(s)$ comme le produit des facteurs $(s+(1+r)/\lambda_0)$
lorsque $r$ parcourt les valuations de ces \'el\'ements $wyf'_y$.

\stl

Pour un multiple plus fin, nous ne consid\'ererons que les poids 
$\rho(w)$ avec $w \in {\cal O} \cap {\cal F}$, 
et $\mbox{val}(w) < 0$. 

\subsubsection{Le polyn\^ome $M_{gros}^C(s)$}

Constatons que pour tout $w \in {\cal F}_{\pi}$, 
 l'\'el\'ement $swf^{s-1}$  appartient \`a  
${\cal D}[s] {\cal F}_0 f^s$. Cela 
r\'esulte de la d\'ecomposition : $w = uf'_x + vf'_y$  
\'etablie au  corollaire \ref{cororesolw}, et de l'identit\'e~:
$$ swf^{s-1} = -(u'_x + v'_y)f^s + \frac{\partial}{\partial x} uf^s + 
\frac{\partial}{\partial y} vf^s\ .$$

\stl

Pour avoir un multiple grossier $M^C_{gros}(s)$ en classe C, 
nous prendrons le produit des polynômes $M_r(s-1)$ pour
tous les poids $r$ parcourant $(1/d){\bf N}$ avec 
$\delta < r < \pi$. Le fait qu'il
n'apparaisse en classe C que des \'el\'ements de poids 
strictement sup\'erieur \`a $\delta$ est une cons\'equence 
du lemme \ref{lemB}. 

\subsubsection{Un multiple de la $b$-fonction}

Notons $\Lambda$ l'ensemble des valeurs propres non nulles de
la matrice magique associ\'ee \`a $f$. D'apr\`es la proposition
\ref{grnabla}, le polyn\^ome minimal de l'action de $\mbox{gr}\,\nabla$
sur $\mbox{gr}_r\,\overline{\cal F}$  est alors~:
$M_r(s)=\prod_{\lambda\in \Lambda}(s+1+(1+r)/\lambda)$. 
Nous rappelons encore la d\'efinition des multiples grossiers
en chaque classe~:

\begin{eqnarray*}
M^A_{gros}(s) &=& \prod_{r\in[0,\pi[\cap(1/d){\bf N}}M_r(s)\ , \\
M^C_{gros}(s) &=& \prod_{r\in]\delta,\pi[\cap(1/d){\bf N}}M_r(s-1)\ , \\
\end{eqnarray*}
$$M^B_{gros}(s)\ =\ \prod_{0\leq l\leq n-2,\ 0\leq k+l\delta<{\sup}\{m_i\}}
(s+\frac{k+1}{n\delta}+\frac{l+1}{n}) $$
o\`u $d$ est le p.p.c.m des degr\'es des facteurs
irr\'eductibles de $f$, et $k,l$ sont des entiers naturels. 

\begin{theo} \label{theomultopo}
Soit $f\in{\bf C}\{x\}[y]$ un polynôme distingu\'e, r\'eduit, 
de degr\'e et de multiplicit\'e $n$.
Le polynôme $(s+1)M_{gros}^A(s)M_{gros}^B(s)M_{gros}^C(s)$ est un multiple du polynôme de Bernstein de $f$
ne d\'ependant que du type topologique du germe de courbe plane d\'efini par $f$.
\end{theo}

\begin{demo}
\`A partir de l'\'el\'ement $ f^s $ de la classe B, nous sommes en mesure 
de faire monter les poids en
utilisant les lemmes \ref{lemB} ou \ref{lemC}. Lors de leur utilisation,
nous avons besoin de comparer les poids de deux \'el\'ements de la 
forme $w'f^s$ et $sw'{'}f^{s-1}$; pour que ça colle, nous devons poser:
$$\varrho (w'f^s , sw'{'}f^{s-1}) \; = \; \mbox{inf} \{ \, n\delta \rho(w') \; , \; \mbox{val}(w'{'}) \,\} $$
et faire monter le poids $\varrho$ des termes rencontr\'es. 
Nous pouvons remarquer  que s'il
y a \'egalit\'e : $ r = n\delta \rho(w') = \mbox{val}(w'{'}) $, 
et si le facteur $(s + \rho(w') + 1/n\delta + 1/n)$
divise $M_r(s-1)$, la multiplication par $M_r(s-1)$ suffit 
pour faire monter le poids des deux \'el\'ements. Nous obtenons alors :
$$M_{gros}^C(s)M^B_{gros}(s)f^s\in{\cal D}[s]{\cal F}_0f^s +
{\cal D}[s]f'_yf^s$$
et on conclut \`a l'aide de l'identit\'e (\ref{eqgrosA}).
\end{demo}

\begin{rema}{\em
   Nous rappelons que M. Saito montre dans \cite{Sa} que
les racines du polyn\^ome de Bernstein d'un germe r\'eduit
$f\in {\cal O}$ sont strictement comprises  entre $-2$ et $0$. 
On peut donc tronquer les polyn\^omes donn\'es pour obtenir 
un multiple plus petit. Nous n'avons pas r\'eussi, pour le 
moment, \`a retrouver explicitement ces in\'egalit\'es par notre 
m\'ethode. }
\end{rema}


\begin{thebibliography}{23}


\bibitem [1]{Mal2}
{\sc Bony J.-M.},  Polyn\^omes de Bernstein et monodromie 
(d'apr\`es B. Malgrange),  S\'eminaire Bourbaki (1974/1975), 
Exp. No. 459, pp. 97--110.

\bibitem [2]{B}
{\sc Brian\c con J.}, \`A propos d'une question de J. Mather,
Universit\'e de Nice-Sophia Antipolis (1972).

\bibitem [3]{BGMM}
{\sc Brian\c con J., Granger M., Maisonobe Ph., Miniconi M.,}
 Algorithme de calcul du polyn\^ome de Bernstein : cas non 
d\'eg\'en\'er\'e, Ann. Inst. Fourier (Grenoble)  39  (1989) 553--610.

\bibitem [4]{K}
{\sc Kashiwara M.,} $B$-functions and holonomic systems,
  Invent. Math.  38  (1976) 33--53.

\bibitem [5]{Mal1}
{\sc Malgrange B.,}  Le polyn\^ome de Bernstein d'une 
singularit\'e isol\'ee, Lecture Notes in Math., Vol. 459, Springer, Berlin, 
1975, pp. 98--199.

\bibitem [6]{Sa}
{\sc Saito M.}, On microlocal $b$-function,  Bull. Soc. math. France  122  
(1994) 163--184.

\bibitem [7]{Var}
{\sc Varchenko A. N.}, Gauss-Manin connection of isolated 
singular point and Bernstein polynomial,  Bull. Sci. Math.  104  
(1980) 205--223.

\bibitem [8]{Ya}
{\sc Yano T.}, On the theory of $b$-functions,  
Publ. Res. Inst. Math. Sci.  14  (1978) 111--202.

\end{thebibliography}
\end{document}